\documentclass{svjour3}                     
\smartqed  
\usepackage{color}
\usepackage{amsmath,amsfonts,amssymb,bm,mathtools}
\usepackage{graphicx,psfrag,subfigure}
\usepackage{enumitem}
\usepackage{epstopdf}
\usepackage{latexsym}
\usepackage{mathrsfs}
\usepackage{empheq}
\usepackage{algorithm}
\usepackage{caption}

\usepackage{fullpage}

\allowdisplaybreaks

\newtheorem{thm}{Theorem}

\newtheorem{rem}{Remark}

\numberwithin{equation}{section}

\begin{document}

\title{A sturcture-preserving, upwind-SAV scheme for the degenerate Cahn--Hilliard equation with applications to simulating surface diffusion
}

\titlerunning{Upwind-SAV scheme for the degenerate Cahn--Hilliard equation}

\author{Qiong-Ao Huang \and Wei Jiang \and Jerry Zhijian Yang \and Cheng Yuan}


\institute{Qiong-Ao Huang\at
School of Mathematics and Statistics, Henan University, Kaifeng 475004, China\\
Center for Applied Mathematics of Henan Province, Henan University, Zhengzhou 450046, China\\
\email{huangqiongao@henu.edu.cn}
\and
Wei Jiang\and Jerry Zhijian Yang\at
School of Mathematics and Statistics, Wuhan University, Wuhan 430072, China\\
Hubei Key Laboratory of Computational Science, Wuhan University, Wuhan 430072, China\\
\email{jiangwei1007@whu.edu.cn; zjyang.math@whu.edu.cn}
\and
Cheng Yuan\at
School of Mathematics and Statistics, Wuhan University, Wuhan 430072, China\\
\email{yuancheng@whu.edu.cn}
}

\date{Received: date / Accepted: date}

\maketitle

\begin{abstract}
This paper establishes a structure-preserving numerical scheme for the Cahn--Hilliard equation with degenerate mobility. First, by applying a finite volume method with upwind numerical fluxes to
the degenerate Cahn--Hilliard equation rewritten by the scalar auxiliary variable (SAV) approach,
we creatively obtain an unconditionally bound-preserving, energy-stable and fully-discrete scheme,
which, for the first time, addresses the boundedness of the classical SAV approach under $H^{-1}$-gradient flow. Then, a dimensional-splitting technique is introduced in high-dimensional cases,
which greatly reduces the computational complexity while preserves original structural properties.
Numerical experiments are presented to verify the bound-preserving and energy-stable properties of the proposed scheme. Finally, by applying the proposed structure-preserving scheme,
we numerically demonstrate that surface diffusion can be approximated by the Cahn--Hilliard equation with degenerate mobility and Flory--Huggins potential when the absolute temperature is sufficiently low, which agrees well with the theoretical result by using formal asymptotic analysis.


\end{abstract}

\section{Introduction}
\label{secint}

The famous Cahn--Hilliard equation was originally established to model phase separation and coarsening processes in binary alloys~\cite{Cahn58},
and is now widely used in many scientific and engineering fields such as 
solid-state dewetting~\cite{Jiang12,Huang19b}, image inpainting~\cite{Bertozzi07,Burger09},
polymer blends~\cite{Brown92,Muller97} and tumor growth~\cite{Garcke18,Ipocoana22},
all of which are built on the following total free energy with respect to the conserved order parameter (i.e., phase variable) $\phi(\bm{x},t)$,
\begin{equation}
\mathcal{E}[\phi(\bm{x},t)]=\int_{\Omega}\left(\frac{\varepsilon^{2}}{2}|\nabla\phi|^{2}+F(\phi)\right)\mathrm{d}\bm{x},\quad
(\bm{x},t)\in\Omega\times[0,T],\label{eq1}
\end{equation}
where $\Omega\subset\mathbb{R}^{d}$ is an open and bounded domain with the boundary $\partial\Omega$,
$\varepsilon>0$ denotes the thickness of the interface between the two phases and
$F(\phi)$ is the Helmholtz free energy density of the system. A typical thermodynamically relevant form of $F(\phi)$ is the so-called logarithmic Flory--Huggins potential as follows~\cite{Cahn58}
\begin{equation}
F_{log}(\phi)=\frac{\theta}{2}\big[(1+\phi)\ln(1+\phi)+(1-\phi)\ln(1-\phi)\big]+\frac{\theta_{c}}{2}(1-\phi^{2}),\quad 0<\theta<\theta_{c},\label{eq2}
\end{equation}
where $\theta$ and $\theta_{c}$ are the absolute and critical temperatures, respectively.
Furthermore, it is easy to check that \eqref{eq2} has a double-well structure with two minima at $\pm\beta_{_{\theta,\theta_{c}}}\in(-1,1)$ and approach $\pm1$ as $\theta/\theta_{c}\to 0$.
Due to the singularity of the logarithmic function in \eqref{eq2},
a simplified polynomial version is used more frequently, namely, 
\begin{equation}
F_{pol}(\phi)=\frac{1}{4}(1-\phi^{2})^{2}.\label{eq3}
\end{equation}
Similarly, \eqref{eq3} is also a double-well potential with two minima at $\pm 1$.

By applying the conserved $H^{-1}$-gradient flow to the energy functional \eqref{eq1},
we can obtain the following desired Cahn--Hilliard equation
\begin{equation}
\frac{\partial\phi}{\partial t}=-\nabla\cdot\bm{J},\quad
\bm{J}=-M(\phi)\nabla\mu,\quad
\mu=-\varepsilon^{2}\Delta\phi+F'(\phi),\quad
 \text{in}~\Omega\times(0,T],\label{eq4}
\end{equation}
equipped with the Neumann and no-flux boundary conditions:
\begin{equation}
\nabla\phi\cdot\bm{n}=0,\quad\bm{J}\cdot\bm{n}=0,\quad \text{on}~\partial\Omega\times(0,T],\label{eq5}
\end{equation}
where $\bm{J}$ is the mass flux, $M(\phi)\ge 0$ is the diffusion mobility, $\mu$ is the chemical potential
and $\bm{n}$ is the outwards pointing normal unit vector onto $\partial\Omega$.

Equipping different types of diffusion mobility $M(\phi)$ into the Cahn--Hilliard equation \eqref{eq4}-\eqref{eq5},
although not changing the energy landscape, can have a significant effect on the kinetic process of $\phi$.
For constant mobility, Pego~\cite{Pego89} and Alikakos {\it et al.}~\cite{Alikakos94} showed by matched asymptotic analysis that
the sharp-interface limit of the Cahn--Hilliard equation is the Mullins--Sekerka problem \cite{Mullins63},
of which the main feature is that a bulk diffusion will appear along with the surface diffusion.
On the other hand, the mobility function degenerate near the minima of potential is also commonly used, which is typically phase-dependent and takes the form of~\cite{Elliott96,Dai16,Pesce21}
\begin{equation}
M(\phi)=(1-\phi^{2})^{k},\quad k=1,2,3,\cdots.\label{eq6}
\end{equation}
Due to the degeneracy of the above mobility, bulk diffusion will be suppressed
and the kinetics are dominated by surface diffusion along the interface of the two phases~\cite{Cahn94}.
In particular, Lee {\it et al.}~\cite{Lee16} showed by matched asymptotic analysis that the sharp-interface limit of
the Cahn--Hilliard equation with polynomial potential~\eqref{eq3} and degenerate mobility~\eqref{eq6} at $k=2$ is surface diffusion,
at least to leading order.
However, the numerical solution owns some oscillations when pure state is reached,
which causes it to be outside the physical range $[-1,1]$,
and the smaller parts may still be absorbed by larger ones in same phase~\cite{Pesce21,Bretin22}.
Pesce \cite{Pesce21} pointed out that the above undesirable results may be a numerical artifact,
and recommended a high-degeneracy mobility in numerical simulations, at least $k=4$ in \eqref{eq6}.
Moreover, by combining the logarithmic potential~\eqref{eq2} at $\theta_{c}=1$ with mobility~\eqref{eq6} at $k=1$,
Cahn {\it et al.}~\cite{Cahn96} proved that the sharp-interface limit is surface diffusion
when $\theta=\mathcal{O}(\varepsilon^{\alpha}),\,\alpha>0$.
Unfortunately, the singularity of logarithmic function makes it challenging to construct a structure-preserving scheme for this problem.

From a continuous point of view, the Cahn--Hilliard equation with degenerate mobility possesses several important properties.
Firstly, although the maximum bound principle similar to the Allen--Cahn equation has not been established for the degenerate Cahn--Hilliard equation~\cite{Du21,Elliott96,Dai16},
it is still necessary to impose \textbf{boundedness} ($|\cdot|\leq1$) on the phase variable $\phi$,
especially for the logarithmic potential.
Secondly, the definition of conserved $H^{-1}$-gradient flow indicates the following \textbf{law of mass conservation},
\begin{equation}
\frac{\mathrm{d}m(t)}{\mathrm{d}t}=\frac{\mathrm{d}}{\mathrm{d}t}\int_{\Omega}\phi\,\mathrm{d}\bm{x}=\int_{\Omega}\phi_{t}\,\mathrm{d}\bm{x}
 =-\int_{\Omega}\nabla\cdot\bm{J}\,\mathrm{d}\bm{x}=-\int_{\partial\Omega}\bm{J}\cdot\bm{n}\,\mathrm{d}s=0.\label{mass}
\end{equation}
Finally, the gradient flow naturally satisfies the \textbf{law of energy dissipation}:
\begin{align}
\frac{\mathrm{d}\mathcal{E}(t)}{\mathrm{d}t}
&=\int_{\Omega}\left(\varepsilon^{2}\nabla\phi\cdot\nabla\phi_{t}+F'(\phi)\phi_{t}\right)\mathrm{d}\bm{x}
=\int_{\Omega}\left(-\varepsilon^{2}\Delta\phi+F'(\phi)\right)\phi_{t}\,\mathrm{d}\bm{x}\nonumber\\
&=\int_{\Omega}\mu\nabla\cdot(M(\phi)\nabla\mu)\,\mathrm{d}\bm{x}=-\int_{\Omega}M(\phi)|\nabla\mu|^{2}\,\mathrm{d}\bm{x}\le 0.\label{ener}
\end{align}
In order to avoid non-physical effects in the simulations over a long time,
it is highly desirable to design a structure-preserving scheme at fully-discrete level.
Fortunately, while the mass conservation can be achieved by many schemes, several widely used approaches aimed at the unconditional energy stability have also been developed in recent years,
such as the convex splitting approach~\cite{Eyre98,Backofen19}, the exponential time differencing (ETD) approach~\cite{Cox02,Fu22},
the invariant energy quadratization (IEQ) approach~\cite{Yang19,Yang17,Huang19b},
the scalar auxiliary variable (SAV) approach~\cite{Shen18,Shen19,Huang19a} and its variants~\cite{Liu20,Cheng20,JiangM22}
(collectively referred to as classical SAV approach), and the new SAV approach~\cite{HuangF20,HuangF22}.
As for the most challenging issue on the bound/positivity preserving, some explorations, but not limited to, are summarized as follows:

\begin{itemize}
\item Finite element approach~\cite{Barrett99} use a practical finite element approximation scheme with an intentionally designed
        variational inequality, the author embedded the boundedness property naturally into the solution
        for Cahn--Hilliard equation. This approximation is proved to be well-posed and stable,
        while the law of energy dissipation for the scheme still needs to be clarified.
\item Function transform approach~\cite{Jungel01} achieves bound/positivity preserving via function transform
        which usually results in a more complicated transformed equation, and likewise fails to address energy stability and mass conservation.
\item Cut-off approach~\cite{Lu13,Li20} artificially cut off values outside the desired range. 
        The main advantage of this approach is that it can achieve arbitrarily high-order accuracy in time for some situations
        (e.g., for Allen--Cahn equation~\cite{Li20}),
        while the energy stability and mass conservation are usually difficult to guarantee.
\item Implicit-explicit approach~\cite{Tang16,Liao20} adopts implicit-explicit discretization and central difference in time and space respectively, leading to the negative diagonally dominant property of the discrete matrix of the Laplace operator under appropriate boundary conditions. Especially, this property is crucial for establishing the bound-preserving scheme for Allen--Cahn type equations.
\item ETD approach~\cite{Du21,LiJ21} comes from the Duhamel principle with the nonlinear terms
        approximated by polynomial interpolations in time, followed by the exact temporal integration.
        Nowadays, an unconditionally bound-preserving scheme has been achieved for Allen--Cahn type equations,
        which also benefits from the good properties of the Laplace operator in the discrete sense
        obtained by central difference or lumped-mass finite element methods.
\item Convex splitting method~\cite{Chen19,Dong19} treats
        the contractive and expansive parts of the potential as implicit and explicit, respectively.
        The main advantage is that it can simultaneously ensure energy stability, mass conservation and boundedness,
        but it is usually useless for degenerate mobility and potential where convex-concave decomposition is inaccessible.
\item Upwind approach~\cite{Bessemoulin12,AcostaSoba22} applies upwind scheme to deal with the flux term,
        leading to a scheme suitable for degenerate PDEs and can guarantee mass conservation, but usually has first-order accuracy in time. Recently, this method was further combined with the convex splitting strategy to achieve the
        unconditional energy stability~\cite{Bailo21,AcostaSoba22}.
\item Lagrange multiplier approach~\cite{Cheng22a,Cheng22b} introduces a space-time Lagrange multiplier
        and uses the Karush--Kuhn--Tucker conditions as a constraint to achieve bound/positivity preserving,
        which, however, does not guarantee the energy stability, although the mass conservation can be ensured.
\end{itemize}
In order to obtain a structure-preserving method, several works on the combination of SAV method and above explorations have been studied, including
\begin{itemize}
  \item New SAV approach with function transform strategy~\cite{HuangF22} achieve weakened energy stability, mass conservation and arbitrarily high-order accuracy in time,
        while can only be applied to homogeneous boundary conditions.
  \item Original SAV approach with cut-off strategy~\cite{Yang22} constructs unconditionally bound-preserving, energy dissipation scheme with arbitrarily high-order accuracy in time for the Allen--Cahn type equations. However, some unphysical increase in the original energy may exist due to the modification in discrete energy.
  \item Exponential SAV approach with stabilized implicit-explicit strategy~\cite{Ju22a} propose a first-order scheme unconditionally maintains boundedness and energy dissipation for the Allen--Cahn type equations, while the boundedness of the second-order one is constrained by the time step size.
  \item Generalized SAV approach with stabilized ETD strategy~\cite{Ju22b}:
        with the appropriate stabilization terms, the proposed SAV-ETD schemes are unconditionally bound-preserving and energy dissipation for the Allen--Cahn type equations.
\end{itemize}
Nevertheless, most of above bound-preserving schemes established for Allen--Cahn type equations
are difficult to extend to Cahn--Hilliard equations, since the negative biharmonic operator involved in Cahn--Hilliard equation is not negative diagonally dominant of Laplace operator as in Allen--Cahn equation under the same spatial discretization~\cite{Du21,Tang16,Liao20}.
In other words, among the above approaches, only the upwind-convex splitting method given in \cite{Bailo21} can supply us with
an original structure-preserving scheme for the Cahn--Hilliard equation~\eqref{eq4}-\eqref{eq5} with degenerate mobility~\eqref{eq6}
and logarithmic potential~\eqref{eq2} or polynomial potential~\eqref{eq3}.
However, despite the many advantages of the convex splitting approach,
several obvious shortcomings including the difficulty in applying to anisotropic potential and constructing high-order scheme in the time direction, requires us to further develop an effective and simple to be generalized structure-preserving scheme.

Inspired by the fact that the SAV-like approaches works well with various potential,
and one of the Lagrange multiplier-type SAV approach can guarantee the original energy stability,
we choose to replace the convex splitting with the Lagrange multiplier-type SAV in the upwind-convex splitting approach,
resulting in a structure-preserving scheme for the Cahn--Hilliard equation with degenerate mobility,
and further numerically verify whether it can capture the main features of surface diffusion.
In high dimensional case, we introduce the dimensional-splitting technique~\cite{Bailo21} into the spatial discretization,
which greatly reduces the computational effort while maintaining the original structural properties.
To sum up, our main contributions are as follows:
\begin{itemize}
  \item By introducing a upwind scheme in the discretization of classical SAV approach, the boundedness of Cahn--Hilliard equation is additionally realized for the first time.
  \item The dimensional-splitting technique is applied to SAV approach for the first time, which greatly improves the computational efficiency.
  \item We numerically verify the surface diffusion can be modeled by the degenerate Cahn--Hilliard equation with Flory--Huggins potential at low temperature.
\end{itemize}

The rest of this paper is organized as follows.
In Section $2$, the upwind strategy is applied to the Lagrange multiplier-type SAV approach, leading to an unconditionally structure-preserving scheme for the Cahn--Hilliard equation with degenerate mobility. In Section $3$, the dimensional-splitting technique is applied to the situation of high-dimensional space to reduce the computational effort caused by space discretization while ensuring the original structural properties.
In Section $4$, ample examples will be provided to show the boundedness, mass conservation and energy stability of the proposed scheme,
and the claim that the sharp-interface limit of the Cahn--Hilliard equation
with degenerate mobility and logarithmic potential at low temperature is surface diffusion will be verified numerically.
Finally, some conclusions will be given in Section $5$.

\section{Upwind-SAV approach}
\label{sec2}

To begin with, we apply the Lagrange multiplier-type SAV approach~\cite{Cheng20,Huang22MCS}
to the degenerate Cahn--Hilliard equation~\eqref{eq4}-\eqref{eq5}:
\begin{equation}
\left\{\begin{array}{l}
\frac{\partial\phi}{\partial t}=-\nabla\cdot\bm{J},\quad
\bm{J}=-M(\phi)\nabla\mu,\quad
\mu=-\varepsilon^{2}\Delta\phi+\xi(t)F'(\phi),\\[2mm]
\frac{\mathrm{d}}{\mathrm{d}t}\int_{\Omega}F(\phi) \mathrm{d}\bm{x}=\xi(t)\int_{\Omega}F'(\phi) \frac{\partial\phi}{\partial t}\mathrm{d}\bm{x},
\end{array}
\right.\quad \text{in}~\Omega\times(0,T],
\label{eqn1}
\end{equation}
subject to the following Neumann and no-flux boundary conditions
\begin{equation}
\nabla\phi\cdot\bm{n}=0,\quad\bm{J}\cdot\bm{n}=0,\quad \text{on}~\partial\Omega\times(0,T],\label{eqn2}
\end{equation}
where $\xi(t)$ is the newly introduced  scalar auxiliary variable whose role is the usual Lagrange multiplier
and $M(\phi)\ge 0$ is the diffusion mobility defined in \eqref{eq6}.
If the initial condition of $\xi(t)$ is taken as $\xi(0)=1$,
it is easy to check that the above rewritten system is equivalent to the original one, i.e., $\xi(t)\equiv 1$.
Therefore, the above equivalent system~\eqref{eqn1}-\eqref{eqn2} also obeys the properties of mass conservation
and energy dissipation at the PDE level, while the boundedness of $\phi$ should also be an essential requirements to be satisfied.

Secondly, to facilitate subsequent numerical discretization, we can define
\begin{equation}
\chi^{+}=\max\{\chi,0\},\quad\chi^{-}=\min\{\chi,0\},\label{eqn3}
\end{equation}
and express the diffusion mobility given in \eqref{eq6} as
\begin{equation}
M(\chi_{1},\chi_{2})=\big[(1+\chi_{1})^{+}(1-\chi_{2})^{+}\big]^{k},\quad k=1,2,3,\cdots.\label{eqn4}
\end{equation}
In fact, due to the boundedness of phase variable, the above expression style is not substantially different from the original one.

Lastly, the finite volume method with upwind numerical fluxes is considered for the discretization of \eqref{eqn1}-\eqref{eqn2}.
Starting with the one dimensional case,
we divide the computational domain $\overline{\Omega}\,(\triangleq\Omega\cup\partial\Omega)$ into $N_{x}$ cells $C_{i}=[x_{i-\frac{1}{2}},x_{i+\frac{1}{2}}],\,i=1,2,\cdots,N_{x}$,
all with uniform size $\Delta x$, so that the centre of each cell is $x_{i}=(i-\frac{1}{2})\Delta x+x_{\frac{1}{2}}$.
By defining the cell average $\phi_{i}$ on $C_{i}$ as
\begin{equation}
\phi_{i}(t)=\frac{1}{\Delta x}\int_{C_{i}}\phi(x,t)\mathrm{d}x\label{eqn5}
\end{equation}
and using the backward Euler formula with finite volume method in temporal and spatial discretization respectively,
the system~\eqref{eqn1} can be approximated as
\begin{align}
&\phi_{i}^{n+1}-\phi_{i}^{n}=-\frac{\Delta t}{\Delta x}\left(J_{i+\frac{1}{2}}^{n+1}-J_{i-\frac{1}{2}}^{n+1}\right),\label{eqn6}\\[2mm]
&J_{i+\frac{1}{2}}^{n+1}=\left(V_{i+\frac{1}{2}}^{n+1}\right)^{+}M(\phi_{i}^{n+1},\phi_{i+1}^{n+1})
                        +\left(V_{i+\frac{1}{2}}^{n+1}\right)^{-}M(\phi_{i+1}^{n+1},\phi_{i}^{n+1}),\label{eqn7}\\[2mm]
&V_{i+\frac{1}{2}}^{n+1}=-\frac{1}{\Delta x}\left(\mu_{i+1}^{n+1}-\mu_{i}^{n+1}\right),\label{eqn8}\\[2mm]
&\mu_{i}^{n+1}=-\varepsilon^{2}(\Delta\phi)_{i}^{n+1}+\xi^{n+1}F'(\phi_{i}^{n+1}),\label{eqn9}\\[2mm]
&\sum_{i=1}^{N_{x}}\left(F(\phi_{i}^{n+1})-F(\phi_{i}^{n})\right)
  =\xi^{n+1}\sum_{i=1}^{N_{x}}F'(\phi_{i}^{n+1})(\phi_{i}^{n+1}-\phi_{i}^{n}),\label{eqn10}
\end{align}
where in \eqref{eqn7} the key idea of the upwind approach~\cite{Bailo21,Gottlieb01,Yan02} has been used.
Here $\Delta t>0$ denotes the time step size, $\chi_{i}^{n}$ is the numerical approximation of $\chi_{i}(t)$
($\chi$ stands for $\phi, J, V, \mu$, $\xi$)
at time $t=t^{n}\triangleq\Delta t\cdot n$ for $n=0,1,\cdots,N$ with $T=N\Delta t$,
and the Laplacian term $(\Delta\phi)_{i}^{n+1}$ is discretized by the following central difference formula
\begin{align}
(\Delta\phi)_{i}^{n+1}&=\frac{(\nabla\phi)_{i+\frac{1}{2}}^{n+1}-(\nabla\phi)_{i-\frac{1}{2}}^{n+1}}{\Delta x}
=\frac{1}{\Delta x}\left(\frac{\phi_{i+1}^{n+1}-\phi_{i}^{n+1}}{\Delta x}-\frac{\phi_{i}^{n+1}-\phi_{i-1}^{n+1}}{\Delta x}\right).\label{eqn11}
\end{align}
Moreover, the Neumann and no-flux boundary conditions are implemented as
\begin{equation}
(\nabla\phi)_{_\frac{1}{2}}^{n+1}=0,\quad (\nabla\phi)_{_{N_{x}+\frac{1}{2}}}^{n+1}=0,\quad
J_{_\frac{1}{2}}^{n+1}=0,\quad J_{_{N_{x}+\frac{1}{2}}}^{n+1}=0,\label{eqn12}
\end{equation}
resulting in that the Laplacian term at the boundaries can be approximated as
\begin{equation}
(\Delta\phi)_{1}^{n+1}=\frac{\phi_{2}^{n+1}-\phi_{1}^{n+1}}{\Delta x^{2}},\quad
(\Delta\phi)_{_{N_{x}}}^{n+1}=\frac{-\phi_{_{N_{x}}}^{n+1}+\phi_{_{N_{x}-1}}^{n+1}}{\Delta x^{2}}.\label{eqn13}
\end{equation}

Now we show the fully-discrete scheme~\eqref{eqn6}-\eqref{eqn13} is structure-preserving.

\begin{thm}\label{thm1}(Boundedness)
The fully-discrete scheme \eqref{eqn6}-\eqref{eqn13} ensure the boundedness of the phase average $\phi_{i}$,
i.e. for $\forall i$, if $|\phi_{i}^{n}|<1$, then $|\phi_{i}^{n+1}|<1$.
\end{thm}

\begin{proof}
First, for $\forall i$, $|\phi_{i}^{n}|<1$ leads to $\phi_{i}^{n+1}<1$. Otherwise,
suppose there is a group of contiguous cells $\{\phi_{j}^{n+1},\phi_{j+1}^{n+1},\cdots,\phi_{k}^{n+1}\}$ such that $\phi_{i}^{n+1}\ge 1$,
then sum the two ends of \eqref{eqn6} over the cells, resulting in
\begin{align}
0<&\,\frac{\Delta x}{\Delta t}\sum_{i=j}^{k}(\phi_{i}^{n+1}-\phi_{i}^{n})
=-\sum_{i=j}^{k}\left(J_{i+\frac{1}{2}}^{n+1}-J_{i-\frac{1}{2}}^{n+1}\right)
=J_{j-\frac{1}{2}}^{n+1}-J_{k+\frac{1}{2}}^{n+1}\nonumber\\
=&\,\left(V_{j-\frac{1}{2}}^{n+1}\right)^{+}M(\phi_{j-1}^{n+1},\phi_{j}^{n+1})
 +\left(V_{j-\frac{1}{2}}^{n+1}\right)^{-}M(\phi_{j}^{n+1},\phi_{j-1}^{n+1})\nonumber\\
&\,-\left(V_{k+\frac{1}{2}}^{n+1}\right)^{+}M(\phi_{k}^{n+1},\phi_{k+1}^{n+1})
 -\left(V_{k+\frac{1}{2}}^{n+1}\right)^{-}M(\phi_{k+1}^{n+1},\phi_{k}^{n+1}).\label{thb1}
\end{align}
On the other hand, since
\begin{equation*}
M(\phi_{j-1}^{n+1},\phi_{j}^{n+1})=M(\phi_{k+1}^{n+1},\phi_{k}^{n+1})=0,~M(\phi_{j}^{n+1},\phi_{j-1}^{n+1})>0,~
M(\phi_{k}^{n+1},\phi_{k+1}^{n+1})>0, \label{thb2}
\end{equation*}
the right end of \eqref{thb1} should be non-positive. Therefore, there must be $\phi_{i}^{n+1}<1$.
We can prove $\phi_{i}^{n+1}>-1$ following the same procedure.
\end{proof}

\begin{rem}\label{rem1}
If the phase-dependent mobility function degenerates at the extreme point of the logarithmic Flory--Huggins potential~\eqref{eq2},
i.e., $M(\phi)=(\beta_{_{\theta,\theta_{c}}}^{2}-\phi^{2})^{k}$,
it can be similarly obtained from the above Theorem~\ref{thm1} that for any $\|\phi_{i}^{n}\|_{\infty}\le \beta_{_{\theta,\theta_{c}}}$,
the fully-discrete scheme~\eqref{eqn6}-\eqref{eqn13} guarantees $\|\phi_{i}^{n+1}\|_{\infty}\le \beta_{_{\theta,\theta_{c}}}$.
\end{rem}

\begin{thm}\label{thm2} (Mass conservation)
The fully-discrete scheme \eqref{eqn6}-\eqref{eqn13} ensures that the total mass is conserved during the evolution, i.e.
\begin{equation}
\sum_{i=1}^{N_{x}}\phi_{i}^{n+1}=\sum_{i=1}^{N_{x}}\phi_{i}^{n}=\cdots=\sum_{i=1}^{N_{x}}\phi_{i}^{0}.\label{tha1}
\end{equation}
\end{thm}

\begin{proof}
Sum the two ends of \eqref{eqn6} over all cells $C_{i}$, it obtains
\begin{equation}
\sum_{i=1}^{N_{x}}(\phi_{i}^{n+1}-\phi_{i}^{n})
=-\frac{\Delta t}{\Delta x}\sum_{i=1}^{N_{x}}\left(J_{i+\frac{1}{2}}^{n+1}-J_{i-\frac{1}{2}}^{n+1}\right)
=-\frac{\Delta t}{\Delta x}\left(J_{{N_{x}+\frac{1}{2}}}^{n+1}-J_{\frac{1}{2}}^{n+1}\right)=0,\label{tha2}
\end{equation}
where the last equality is derived from the non-flux boundary conditions \eqref{eqn12}.
\end{proof}

\begin{thm}\label{thm3} (Energy dissipation)
The fully-discrete scheme \eqref{eqn6}-\eqref{eqn13} is unconditionally energy stable,
satisfying the following discrete energy dissipation law:
\begin{equation}
\frac{\mathcal{E}^{n+1}-\mathcal{E}^{n}}{\Delta x}\le -\Delta t
\sum_{i=1}^{N_{x}-1}\min\left\{M(\phi_{i}^{n+1},\phi_{i+1}^{n+1}),M(\phi_{i+1}^{n+1},\phi_{i}^{n+1})\right\}
\left|V_{i+\frac{1}{2}}^{n+1}\right|^{2}\le 0,\label{thc1}
\end{equation}
where
\begin{equation}
\mathcal{E}^{n}=\Delta x\sum_{i=1}^{N_{x}-1}\frac{\varepsilon^{2}}{2}\big|(\nabla\phi)_{i+\frac{1}{2}}^{n}\big|^{2}
+\Delta x\sum_{i=1}^{N_{x}}F(\phi_{i}^{n}).\label{thc2}
\end{equation}
\end{thm}

\begin{proof}
Subtracting the discrete free energy in \eqref{thc2} at subsequent times leads to
\begin{align}
\frac{\mathcal{E}^{n+1}-\mathcal{E}^{n}}{\Delta x}
=&\,\frac{\varepsilon^{2}}{2}\sum_{i=1}^{N_{x}-1}
  \left(\big|(\nabla\phi)_{i+\frac{1}{2}}^{n+1}\big|^{2}-\big|(\nabla\phi)_{i+\frac{1}{2}}^{n}\big|^{2}\right)
 +\sum_{i=1}^{N_{x}}\left(F(\phi_{i}^{n+1})-F(\phi_{i}^{n})\right)\nonumber\\
=&\,\varepsilon^{2}\sum_{i=1}^{N_{x}-1}(\nabla\phi)_{i+\frac{1}{2}}^{n+1}
  \left[(\nabla\phi)_{i+\frac{1}{2}}^{n+1}-(\nabla\phi)_{i+\frac{1}{2}}^{n}\right]\nonumber\\
  &\,-\frac{\varepsilon^{2}}{2}\sum_{i=1}^{N_{x}-1}\big|(\nabla\phi)_{i+\frac{1}{2}}^{n+1}-(\nabla\phi)_{i+\frac{1}{2}}^{n}\big|^{2}
  +\xi^{n+1}\sum_{i=1}^{N_{x}}F'(\phi_{i}^{n+1})(\phi_{i}^{n+1}-\phi_{i}^{n})\nonumber\\
=&\,\varepsilon^{2}\sum_{i=1}^{N_{x}-1}\frac{\phi_{i+1}^{n+1}-\phi_{i}^{n+1}}{\Delta x^{2}}(\phi_{i+1}^{n+1}-\phi_{i+1}^{n})
    -\varepsilon^{2}\sum_{i=1}^{N_{x}-1}\frac{\phi_{i+1}^{n+1}-\phi_{i}^{n+1}}{\Delta x^{2}}(\phi_{i}^{n+1}-\phi_{i}^{n})\nonumber\\
  &\,-\frac{\varepsilon^{2}}{2}\sum_{i=1}^{N_{x}-1}\big|(\nabla\phi)_{i+\frac{1}{2}}^{n+1}-(\nabla\phi)_{i+\frac{1}{2}}^{n}\big|^{2}
     +\xi^{n+1}\sum_{i=1}^{N_{x}}F'(\phi_{i}^{n+1})(\phi_{i}^{n+1}-\phi_{i}^{n})\nonumber\\
=&\,-\varepsilon^{2}\sum_{i=2}^{N_{x}-1}\frac{\phi_{i+1}^{n+1}-2\phi_{i}^{n+1}+\phi_{i-1}^{n+1}}{\Delta x^{2}}(\phi_{i}^{n+1}-\phi_{i}^{n})
  \nonumber\\
  &\,+\varepsilon^{2}\frac{\phi_{_{N_{x}}}^{n+1}-\phi_{_{N_{x}-1}}^{n+1}}{\Delta x^{2}}(\phi_{_{N_{x}}}^{n+1}-\phi_{_{N_{x}}}^{n})
     -\varepsilon^{2}\frac{\phi_{2}^{n+1}-\phi_{1}^{n+1}}{\Delta x^{2}}(\phi_{1}^{n+1}-\phi_{1}^{n})\nonumber\\
  &\,-\frac{\varepsilon^{2}}{2}\sum_{i=1}^{N_{x}-1}\big|(\nabla\phi)_{i+\frac{1}{2}}^{n+1}-(\nabla\phi)_{i+\frac{1}{2}}^{n}\big|^{2}
     +\xi^{n+1}\sum_{i=1}^{N_{x}}F'(\phi_{i}^{n+1})(\phi_{i}^{n+1}-\phi_{i}^{n})\nonumber\\
\le&\,-\sum_{i=1}^{N_{x}}\varepsilon^{2}(\Delta\phi)_{i}^{n+1}(\phi_{i}^{n+1}-\phi_{i}^{n})
  +\xi^{n+1}\sum_{i=1}^{N_{x}}F'(\phi_{i}^{n+1})(\phi_{i}^{n+1}-\phi_{i}^{n})\nonumber\\
=&\,-\sum_{i=1}^{N_{x}}\frac{\Delta t}{\Delta x}\left(J_{i+\frac{1}{2}}^{n+1}-J_{i-\frac{1}{2}}^{n+1}\right)\mu_{i}^{n+1}
=\Delta t\sum_{i=1}^{N_{x}-1}J_{i+\frac{1}{2}}^{n+1}\frac{\mu_{i+1}^{n+1}-\mu_{i}^{n+1}}{\Delta x}\nonumber\\
=&\,-\Delta t\sum_{i=1}^{N_{x}-1}\left[\big(V_{i+\frac{1}{2}}^{n+1}\big)^{+}M(\phi_{i}^{n+1},\phi_{i+1}^{n+1})
                        +\big(V_{i+\frac{1}{2}}^{n+1}\big)^{-}M(\phi_{i+1}^{n+1},\phi_{i}^{n+1})\right]V_{i+\frac{1}{2}}^{n+1}\nonumber\\
\le&\,-\Delta t\sum_{i=1}^{N_{x}-1}\min\left\{M(\phi_{i}^{n+1},\phi_{i+1}^{n+1}),M(\phi_{i+1}^{n+1},\phi_{i}^{n+1})\right\}
     \left|V_{i+\frac{1}{2}}^{n+1}\right|^{2}\le\,0.\label{thc3}
\end{align}
\end{proof}

Theorem~\ref{thm1} shows that the key to establishing the boundedness of numerical solution is to use the upwind idea to deal with the mass flux numerically, and its bound is only determined by the zero point of the degenerate mobility function.
Meanwhile, SAV approach ensure the stability of energy and the conservation of mass.

On the other hand, although the above scheme \eqref{eqn6}-\eqref{eqn13} can be directly extended to high dimensional case, we would introduce the dimensional-splitting technique~\cite{Bailo21} in order to save computational effort.

\section{Dimensional-splitting technique}
\label{sec3}

We begin with dividing the computational domain $\overline{\Omega}$ uniformly into $N_{x}\times N_{y}$ cells
$C_{i,j}=[x_{i-\frac{1}{2}},x_{i+\frac{1}{2}}]\times[y_{j-\frac{1}{2}},y_{j+\frac{1}{2}}]$,
with the spatial steps $\Delta x$ and $\Delta y$. In each cell $C_{i,j}$, the cell average $\phi_{i,j}$ is defined as
\begin{equation}
\phi_{i,j}(t)=\frac{1}{\Delta x\Delta y}\int_{C_{i,j}}\phi(x,y,t)\mathrm{d}x\mathrm{d}y.\label{eqs1}
\end{equation}

According to the dimensional-splitting technique, at each step $n$ in the outer loop, we first update $\phi_{i,j}$ along $x$-direction for each fixed $y_j$ one by one, then solve the systems along $y$-direction for every fixed $x_i$. More precisely, we can use $\widetilde{\phi}_{i,j}^{n,q}$ to stands for the solution at the $q$-th update in the first inner loop while $q\,(=1,2,\cdots,N_{y})$ denotes the index of the fixed $y_q$ in this loop. Similarly, the solution at the $p$-th update in the second inner loop can be written as  $\widehat{\phi}_{i,j}^{n,p}$ while $p=1,2,\cdots,N_{x}$. Based on these notations, the initial condition of the first inner loop can be expressed as $\widetilde{\phi}_{i,j}^{n,q}|_{q=0}=\phi_{i,j}^{n}$ and the dimensional-splitting technique can be formulated as:

$\mathbf{Step~1.~for~q=1,2,\cdots,N_{y}~do:}$
\begin{align}
&\widetilde{\phi}_{i,j}^{n,q}-\widetilde{\phi}_{i,j}^{n,q-1}=
\left\{\begin{array}{ll}
-\frac{\Delta t}{\Delta x}\left(\widetilde{J}_{i+\frac{1}{2},j}^{n,q}-\widetilde{J}_{i-\frac{1}{2},j}^{n,q}\right),& if~j=q;\\[2mm]
0,& otherwise,
\end{array}\right.\label{eqs2}\\[2mm]
&\widetilde{J}_{i+\frac{1}{2},j}^{n,q}=\left(\widetilde{V}_{i+\frac{1}{2},j}^{n,q}\right)^{+}M(\widetilde{\phi}_{i,j}^{n,q},\widetilde{\phi}_{i+1,j}^{n,q})
                        +\left(\widetilde{V}_{i+\frac{1}{2},j}^{n,q}\right)^{-}M(\widetilde{\phi}_{i+1,j}^{n,q},\widetilde{\phi}_{i,j}^{n,q}),\label{eqs3}\\[2mm]
&\widetilde{V}_{i+\frac{1}{2},j}^{n,q}=-\frac{1}{\Delta x}\left(\widetilde{\mu}_{i+1,j}^{n,q}-\widetilde{\mu}_{i,j}^{n,q}\right),\label{eqs4}\\[2mm]
&\widetilde{\mu}_{i,j}^{n,q}=-\varepsilon^{2}(\Delta\widetilde{\phi})_{i,j}^{n,q}+\widetilde{\xi}^{n,q}F'(\widetilde{\phi}_{i,j}^{n,q}),\label{eqs5}\\[2mm]
&\sum_{i=1}^{N_{x}}\sum_{j=1}^{N_{y}}\left(F(\widetilde{\phi}_{i,j}^{n,q})-F(\widetilde{\phi}_{i,j}^{n,q-1})\right)
  =\widetilde{\xi}^{n,q}\sum_{i=1}^{N_{x}}\sum_{j=1}^{N_{y}}F'(\widetilde{\phi}_{i,j}^{n,q})(\widetilde{\phi}_{i,j}^{n,q}-\widetilde{\phi}_{i,j}^{n,q-1}),\label{eqs6}
\end{align}
where  the Laplacian term $(\Delta\widetilde{\phi})_{i,j}^{n,q}$ at the boundaries are discretized as
\begin{equation}
\left\{\begin{array}{l}
(\Delta\widetilde{\phi})_{1,1}^{n,q}=\frac{\widetilde{\phi}_{2,1}^{n,q}-\widetilde{\phi}_{1,1}^{n,q}}{\Delta x^{2}}
   +\frac{\widetilde{\phi}_{1,2}^{n,q}-\widetilde{\phi}_{1,1}^{n,q}}{\Delta y^{2}},\\[2mm]
(\Delta\widetilde{\phi})_{N_{x},1}^{n,q}=\frac{-\widetilde{\phi}_{N_{x},1}^{n,q}+\widetilde{\phi}_{N_{x}-1,1}^{n,q}}{\Delta x^{2}}
   +\frac{\widetilde{\phi}_{N_{x},2}^{n,q}-\widetilde{\phi}_{N_{x},1}^{n,q}}{\Delta y^{2}},\\[2mm]
(\Delta\widetilde{\phi})_{1,N_{y}}^{n,q}=\frac{\widetilde{\phi}_{2,N_{y}}^{n,q}-\widetilde{\phi}_{1,N_{y}}^{n,q}}{\Delta x^{2}}
   +\frac{-\widetilde{\phi}_{1,N_{y}}^{n,q}+\widetilde{\phi}_{1,N_{y}-1}^{n,q}}{\Delta y^{2}},\\[2mm]
(\Delta\widetilde{\phi})_{N_{x},N_{y}}^{n,q}=\frac{-\widetilde{\phi}_{N_{x},N_{y}}^{n,q}+\widetilde{\phi}_{N_{x}-1,N_{y}}^{n,q}}{\Delta x^{2}}
   +\frac{-\widetilde{\phi}_{N_{x},N_{y}}^{n,q}+\widetilde{\phi}_{N_{x},N_{y}-1}^{n,q}}{\Delta y^{2}},\\[2mm]
(\Delta\widetilde{\phi})_{1,j}^{n,q}=\frac{\widetilde{\phi}_{2,j}^{n,q}-\widetilde{\phi}_{1,j}^{n,q}}{\Delta x^{2}}
   +\frac{\widetilde{\phi}_{1,j+1}^{n,q}-2\widetilde{\phi}_{1,j}^{n,q}+\widetilde{\phi}_{1,j-1}^{n,q}}{\Delta y^{2}},\quad j=2,3,\cdots,N_{y}-1,\\[2mm]
(\Delta\widetilde{\phi})_{N_{x},j}^{n,q}=\frac{-\widetilde{\phi}_{N_{x},j}^{n,q}+\widetilde{\phi}_{N_{x}-1,j}^{n,q}}{\Delta x^{2}}
   +\frac{\widetilde{\phi}_{N_{x},j+1}^{n,q}-2\widetilde{\phi}_{N_{x},j}^{n,q}+\widetilde{\phi}_{N_{x},j-1}^{n,q}}{\Delta y^{2}},\quad j=2,3,\cdots,N_{y}-1,\\[2mm]
(\Delta\widetilde{\phi})_{i,1}^{n,q}=\frac{\widetilde{\phi}_{i+1,1}^{n,q}-2\widetilde{\phi}_{i,1}^{n,q}+\widetilde{\phi}_{i-1,1}^{n,q}}{\Delta x^{2}}
   +\frac{\widetilde{\phi}_{i,2}^{n,q}-\widetilde{\phi}_{i,1}^{n,q}}{\Delta y^{2}},\quad i=2,3,\cdots,N_{x}-1,\\[2mm]
(\Delta\widetilde{\phi})_{i,N_{y}}^{n,q}=\frac{\widetilde{\phi}_{i+1,N_{y}}^{n,q}-2\widetilde{\phi}_{i,N_{y}}^{n,q}+\widetilde{\phi}_{i-1,N_{y}}^{n,q}}{\Delta x^{2}}
   +\frac{-\widetilde{\phi}_{i,N_{y}}^{n,q}+\widetilde{\phi}_{i,N_{y}-1}^{n,q}}{\Delta y^{2}},\quad i=2,3,\cdots,N_{x}-1,
\end{array}
\right.\label{eqs8}
\end{equation}
while in the domain defined by the following central difference formula
\begin{equation}
(\Delta\widetilde{\phi})_{i,j}^{n,q}=\frac{\widetilde{\phi}_{i+1,j}^{n,q}-2\widetilde{\phi}_{i,j}^{n,q}+\widetilde{\phi}_{i-1,j}^{n,q}}{\Delta x^{2}}
+\frac{\widetilde{\phi}_{i,j+1}^{n,q}-2\widetilde{\phi}_{i,j}^{n,q}+\widetilde{\phi}_{i,j-1}^{n,q}}{\Delta y^{2}}.\label{eqs7}
\end{equation}
Moreover, the no-flux boundary conditions for mass flux are implemented as
\begin{equation}
\widetilde{J}_{{\frac{1}{2},j}}^{n,q}=0,\quad \widetilde{J}_{N_{x}+\frac{1}{2},j}^{n,q}=0,\quad j=1,2,\cdots,N_{y}.\label{eqs9}
\end{equation}

Once the above first inner loop is completed, the cell average for each of the $y$-direction is continued.
The initial condition for second inner loop is taken as $\widehat{\phi}_{i,j}^{n,p}|_{p=0}=\widetilde{\phi}_{i,j}^{n,q}|_{q=_{N_{y}}}$
and the scheme for each $y$-direction satisfies:

$\mathbf{Step~2.~for~p=1,2,\cdots,N_{x}~do:}$
\begin{align}
&\widehat{\phi}_{i,j}^{n,p}-\widehat{\phi}_{i,j}^{n,p-1}=
\left\{\begin{array}{ll}
-\frac{\Delta t}{\Delta y}\left(\widehat{J}_{i,j+\frac{1}{2}}^{n,p}-\widehat{J}_{i,j-\frac{1}{2}}^{n,p}\right),& if~i=p;\\[2mm]
0,& otherwise,
\end{array}\right.\label{eqs10}\\[2mm]
&\widehat{J}_{i,j+\frac{1}{2}}^{n,p}=\left(\widehat{V}_{i,j+\frac{1}{2}}^{n,p}\right)^{+}M(\widehat{\phi}_{i,j}^{n,p},\widehat{\phi}_{i,j+1}^{n,p})
                        +\left(\widehat{V}_{i,j+\frac{1}{2}}^{n,p}\right)^{-}M(\widehat{\phi}_{i,j+1}^{n,p},\widehat{\phi}_{i,j}^{n,p}),\label{eqs11}\\[2mm]
&\widehat{V}_{i,j+\frac{1}{2}}^{n,p}=-\frac{1}{\Delta y}\left(\widehat{\mu}_{i,j+1}^{n,p}-\widehat{\mu}_{i,j}^{n,p}\right),\label{eqs12}\\[2mm]
&\widehat{\mu}_{i,j}^{n,p}=-\varepsilon^{2}(\Delta\widehat{\phi})_{i,j}^{n,p}+\widehat{\xi}^{n,p}F'(\widehat{\phi}_{i,j}^{n,p}),\label{eqs13}\\[2mm]
&\sum_{i=1}^{N_{x}}\sum_{j=1}^{N_{y}}\left(F(\widehat{\phi}_{i,j}^{n,p})-F(\widehat{\phi}_{i,j}^{n,p-1})\right)
  =\widehat{\xi}^{n,p}\sum_{i=1}^{N_{x}}\sum_{j=1}^{N_{y}}F'(\widehat{\phi}_{i,j}^{n,p})(\widehat{\phi}_{i,j}^{n,p}-\widehat{\phi}_{i,j}^{n,p-1}),\label{eqs14}
\end{align}
where the Laplacian term $(\Delta\widehat{\phi})_{i,j}^{n,p}$ is discretized by the following formula in the domain
\begin{equation}
(\Delta\widehat{\phi})_{i,j}^{n,p}=\frac{\widehat{\phi}_{i+1,j}^{n,p}-2\widehat{\phi}_{i,j}^{n,p}+\widehat{\phi}_{i-1,j}^{n,p}}{\Delta x^{2}}
+\frac{\widehat{\phi}_{i,j+1}^{n,p}-2\widehat{\phi}_{i,j}^{n,p}+\widehat{\phi}_{i,j-1}^{n,p}}{\Delta y^{2}},\label{eqs15}
\end{equation}
and at the boundaries are discretized as
\begin{equation}
\left\{\begin{array}{ll}
(\Delta\widehat{\phi})_{1,1}^{n,p}=\frac{\widehat{\phi}_{2,1}^{n,p}-\widehat{\phi}_{1,1}^{n,p}}{\Delta x^{2}}
   +\frac{\widehat{\phi}_{1,2}^{n,p}-\widehat{\phi}_{1,1}^{n,p}}{\Delta y^{2}},\\[2mm]
(\Delta\widehat{\phi})_{N_{x},1}^{n,p}=\frac{-\widehat{\phi}_{N_{x},1}^{n,p}+\widehat{\phi}_{N_{x}-1,1}^{n,p}}{\Delta x^{2}}
   +\frac{\widehat{\phi}_{N_{x},2}^{n,p}-\widehat{\phi}_{N_{x},1}^{n,p}}{\Delta y^{2}},\\[2mm]
(\Delta\widehat{\phi})_{1,N_{y}}^{n,p}=\frac{\widehat{\phi}_{2,N_{y}}^{n,p}-\widehat{\phi}_{1,N_{y}}^{n,p}}{\Delta x^{2}}
   +\frac{-\widehat{\phi}_{1,N_{y}}^{n,p}+\widehat{\phi}_{1,N_{y}-1}^{n,p}}{\Delta y^{2}},\\[2mm]
(\Delta\widehat{\phi})_{N_{x},N_{y}}^{n,p}=\frac{-\widehat{\phi}_{N_{x},N_{y}}^{n,p}+\widehat{\phi}_{N_{x}-1,N_{y}}^{n,p}}{\Delta x^{2}}
   +\frac{-\widehat{\phi}_{N_{x},N_{y}}^{n,p}+\widehat{\phi}_{N_{x},N_{y}-1}^{n,p}}{\Delta y^{2}},\\[2mm]
(\Delta\widehat{\phi})_{1,j}^{n,p}=\frac{\widehat{\phi}_{2,j}^{n,p}-\widehat{\phi}_{1,j}^{n,p}}{\Delta x^{2}}
   +\frac{\widehat{\phi}_{1,j+1}^{n,p}-2\widehat{\phi}_{1,j}^{n,p}+\widehat{\phi}_{1,j-1}^{n,p}}{\Delta y^{2}},\quad j=2,3,\cdots,N_{y}-1,\\[2mm]
(\Delta\widehat{\phi})_{N_{x},j}^{n,p}=\frac{-\widehat{\phi}_{N_{x},j}^{n,p}+\widehat{\phi}_{N_{x}-1,j}^{n,p}}{\Delta x^{2}}
   +\frac{\widehat{\phi}_{N_{x},j+1}^{n,p}-2\widehat{\phi}_{N_{x},j}^{n,p}+\widehat{\phi}_{N_{x},j-1}^{n,p}}{\Delta y^{2}},\quad j=2,3,\cdots,N_{y}-1,\\[2mm]
(\Delta\widehat{\phi})_{i,1}^{n,p}=\frac{\widehat{\phi}_{i+1,1}^{n,p}-2\widehat{\phi}_{i,1}^{n,p}+\widehat{\phi}_{i-1,1}^{n,p}}{\Delta x^{2}}
   +\frac{\widehat{\phi}_{i,2}^{n,p}-\widehat{\phi}_{i,1}^{n,p}}{\Delta y^{2}},\quad i=2,3,\cdots,N_{x}-1,\\[2mm]
(\Delta\widehat{\phi})_{i,N_{y}}^{n,p}=\frac{\widehat{\phi}_{i+1,N_{y}}^{n,p}-2\widehat{\phi}_{i,N_{y}}^{n,p}+\widehat{\phi}_{i-1,N_{y}}^{n,p}}{\Delta x^{2}}
   +\frac{-\widehat{\phi}_{i,N_{y}}^{n,p}+\widehat{\phi}_{i,N_{y}-1}^{n,p}}{\Delta y^{2}},\quad i=2,3,\cdots,N_{x}-1.
\end{array}
\right.\label{eqs16}
\end{equation}
Moreover, the no-flux boundary conditions for mass flux are implemented as
\begin{equation}
\widehat{J}_{i,\frac{1}{2}}^{n,p}=0,\quad \widehat{J}_{i,N_{y}+\frac{1}{2}}^{n,p}=0,\quad i=1,2,\cdots,N_{x}.\label{eqs17}
\end{equation}

Finally, when the above inner loops are completed,
the cell average and the Lagrange multiplier at the $(n+1)$-th step can be obtained as follows
\begin{equation}
\phi_{i,j}^{n+1}=\widehat{\phi}_{i,j}^{n,p}|_{p=_{N_{x}}}\quad   \text{and}\quad
\xi^{n+1}=\frac{\sum_{q=1}^{N_{y}}\widetilde{\xi}^{n,q}+\sum_{p=1}^{N_{x}}\widehat{\xi}^{n,p}}{N_{y}+N_{x}}.\label{eqs18}
\end{equation}
Notice here we simply use the average of intermediate value of $\xi$ to define the Lagrange multiplier at $(n+1)$-th step. In addition to being efficiency, a more significant virtue of the dimensional-splitting technique is that the structure-preserving property still holds, which can be rigorously proved as follows.

\begin{thm}\label{thm4} (Boundedness)
The dimensional-splitting scheme \eqref{eqs2}-\eqref{eqs18} can ensure the boundedness of $\phi_{i,j}$, i.e. $\forall i,j$, if $|\phi_{i,j}^{n}|<1$, then $|\phi_{i,j}^{n+1}|<1$.
\end{thm}

\begin{proof}
Since $\left|\widetilde{\phi}_{i,j}^{n,q}|_{q=0}\right|=|\phi_{i,j}^{n}|<1$,
it is easy to prove that $|\widetilde{\phi}_{i,j}^{n,q}|<1$ for $\forall q$
by following the contradiction strategy employed for one-dimension case in Theorem~\ref{thm1}.
Similarly, $|\widehat{\phi}_{i,j}^{n,p}|<1$ for $\forall p$ can be obtained from $\left|\widehat{\phi}_{i,j}^{n,p}|_{p=0}\right|=\left|\widetilde{\phi}_{i,j}^{n,q}|_{q=_{N_{y}}}\right|<1$,
and thus $|\phi_{i,j}^{n+1}|=\left|\widehat{\phi}_{i,j}^{n,p}|_{p=_{N_{x}}}\right|<1$.
\end{proof}

\begin{thm}\label{thm5} (Mass conservation)
The dimensional-splitting scheme \eqref{eqs2}-\eqref{eqs18} ensures that the total mass is conserved during the evolution, i.e.
\begin{equation}
\sum_{i=1}^{N_{x}}\sum_{j=1}^{N_{y}}\phi_{i,j}^{n+1}=\sum_{i=1}^{N_{x}}\sum_{j=1}^{N_{y}}\phi_{i,j}^{n}
=\cdots=\sum_{i=1}^{N_{x}}\sum_{j=1}^{N_{y}}\phi_{i,j}^{0}.\label{thd1}
\end{equation}
\end{thm}

\begin{proof}
Sum the two ends of \eqref{eqs2} and \eqref{eqs10} over all cells $C_{i,j}$, respectively, it obtains
\begin{align}
\sum_{i=1}^{N_{x}}\sum_{j=1}^{N_{y}}(\widetilde{\phi}_{i,j}^{n,q}-\widetilde{\phi}_{i,j}^{n,q-1})
&=-\frac{\Delta t}{\Delta x}\sum_{i=1}^{N_{x}}\left(\widetilde{J}_{i+\frac{1}{2},q}^{n,q}-\widetilde{J}_{i-\frac{1}{2},q}^{n,q}\right)=0,
\quad for~\forall q,\label{thd2}\\
\sum_{j=1}^{N_{y}}\sum_{i=1}^{N_{x}}(\widehat{\phi}_{i,j}^{n,p}-\widehat{\phi}_{i,j}^{n,p-1})
&=-\frac{\Delta t}{\Delta y}\sum_{j=1}^{N_{y}}\left(\widehat{J}_{p,j+\frac{1}{2}}^{n,p}-\widehat{J}_{p,j-\frac{1}{2}}^{n,p}\right)=0,
\quad for~\forall p.\label{thd3}
\end{align}
Therefore,
\begin{align}
\sum_{i=1}^{N_{x}}\sum_{j=1}^{N_{y}}\phi_{i,j}^{n+1}
=&\sum_{i=1}^{N_{x}}\sum_{j=1}^{N_{y}}\widehat{\phi}_{i,j}^{n,p}|_{p=_{N_{x}}}
=\sum_{i=1}^{N_{x}}\sum_{j=1}^{N_{y}}\widehat{\phi}_{i,j}^{n,p}|_{p=_{N_{x}}-1}
=\cdots=\sum_{i=1}^{N_{x}}\sum_{j=1}^{N_{y}}\widehat{\phi}_{i,j}^{n,p}|_{p=0}\nonumber\\
=&\sum_{i=1}^{N_{x}}\sum_{j=1}^{N_{y}}\widetilde{\phi}_{i,j}^{n,q}|_{q=_{N_{y}}}
=\sum_{i=1}^{N_{x}}\sum_{j=1}^{N_{y}}\widetilde{\phi}_{i,j}^{n,q}|_{q=_{N_{y}}-1}
=\cdots=\sum_{i=1}^{N_{x}}\sum_{j=1}^{N_{y}}\widetilde{\phi}_{i,j}^{n,q}|_{q=0}\nonumber\\
=&\sum_{i=1}^{N_{x}}\sum_{i=j}^{N_{y}}\phi_{i,j}^{n}
=\cdots=\sum_{i=1}^{N_{x}}\sum_{j=1}^{N_{y}}\phi_{i,j}^{0},\label{thd4}
\end{align}
where the initial conditions $\widehat{\phi}_{i,j}^{n,p}|_{p=0}=\widetilde{\phi}_{i,j}^{n,q}|_{q=_{N_{y}}}$ and $\widetilde{\phi}_{i,j}^{n,q}|_{q=0}=\phi_{i,j}^{n}$ are applied.
\end{proof}

\begin{thm}\label{thm6} (Energy dissipation)
The dimensional-splitting scheme \eqref{eqs2}-\eqref{eqs18} is unconditionally energy stable,
and satisfies the following discrete energy dissipation law:
\begin{align}
\frac{\mathcal{E}^{n+1}-\mathcal{E}^{n}}{\Delta x\Delta y}
\le&\, -\Delta t\sum_{p=1}^{N_{x}}\sum_{j=1}^{N_{y}-1}
\min\left\{M(\widehat{\phi}_{p,j}^{n,p},\widehat{\phi}_{p,j+1}^{n,p}),M(\widehat{\phi}_{p,j+1}^{n,p},\widehat{\phi}_{p,j}^{n,p})\right\}
   \left|\widehat{V}_{p,j+\frac{1}{2}}^{n,p}\right|^{2}\nonumber\\
&\,-\Delta t\sum_{q=1}^{N_{y}}\sum_{i=1}^{N_{x}-1}
\min\left\{M(\widetilde{\phi}_{i,q}^{n,q},\widetilde{\phi}_{i+1,q}^{n,q}),M(\widetilde{\phi}_{i+1,q}^{n,q},\widetilde{\phi}_{i,q}^{n,q})\right\}
   \left|\widetilde{V}_{i+\frac{1}{2},q}^{n,q}\right|^{2}\leq 0,\label{thf1}
\end{align}
where
\begin{align}
\mathcal{E}^{n}=\,&
\Delta x\Delta y\sum_{i=1}^{N_{x}-1}\sum_{j=1}^{N_{y}}\frac{\varepsilon^{2}}{2}\left(\frac{\phi_{i+1,j}^{n}-\phi_{i,j}^{n}}{\Delta x}\right)^{2}
+\Delta x\Delta y\sum_{i=1}^{N_{x}}\sum_{j=1}^{N_{y}-1}\frac{\varepsilon^{2}}{2}\left(\frac{\phi_{i,j+1}^{n}-\phi_{i,j}^{n}}{\Delta y}\right)^{2}
+\Delta x\Delta y\sum_{i=1}^{N_{x}}\sum_{j=1}^{N_{y}}F(\phi_{i,j}^{n}).\label{thf2}
\end{align}
\end{thm}

\begin{proof}
We first show the energy dissipation in each $x$-direction iteration, namely,
\begin{equation}
\widetilde{\mathcal{E}}^{n,q}-\widetilde{\mathcal{E}}^{n,q-1}
\le -\Delta t\Delta x\Delta y\sum_{i=1}^{N_{x}-1}
\min\Big\{M(\widetilde{\phi}_{i,q}^{n,q},\widetilde{\phi}_{i+1,q}^{n,q}),M(\widetilde{\phi}_{i+1,q}^{n,q},\widetilde{\phi}_{i,q}^{n,q})\Big\}
     \Big|\widetilde{V}_{i+\frac{1}{2},q}^{n,q}\Big|^{2}\le 0,\label{thf3}
\end{equation}
where
\begin{align}
\widetilde{\mathcal{E}}^{n,q}=&\,
\Delta x\Delta y\sum_{i=1}^{N_{x}-1}\sum_{j=1}^{N_{y}}\frac{\varepsilon^{2}}{2}
  \left(\frac{\widetilde{\phi}_{i+1,j}^{n,q}-\widetilde{\phi}_{i,j}^{n,q}}{\Delta x}\right)^{2}
+\Delta x\Delta y\sum_{i=1}^{N_{x}}\sum_{j=1}^{N_{y}-1}\frac{\varepsilon^{2}}{2}
  \left(\frac{\widetilde{\phi}_{i,j+1}^{n,q}-\widetilde{\phi}_{i,j}^{n,q}}{\Delta y}\right)^{2}\nonumber\\
&\,+\Delta x\Delta y\sum_{i=1}^{N_{x}}\sum_{j=1}^{N_{y}}F(\widetilde{\phi}_{i,j}^{n,q}).\label{thf4}
\end{align}
By subtracting the first two terms on the right side of \eqref{thf4} at subsequent times,
\begin{align}
&\,\frac{\varepsilon^{2}}{2}\sum_{i=1}^{N_{x}-1}\sum_{j=1}^{N_{y}}\left[\left(\frac{\widetilde{\phi}_{i+1,j}^{n,q}-\widetilde{\phi}_{i,j}^{n,q}}{\Delta x}\right)^{2}
       -\left(\frac{\widetilde{\phi}_{i+1,j}^{n,q-1}-\widetilde{\phi}_{i,j}^{n,q-1}}{\Delta x}\right)^{2}\right]\nonumber\\
=&\,-\varepsilon^{2}\sum_{i=2}^{N_{x}-1}\sum_{j=1}^{N_{y}}
   \left(\frac{\widetilde{\phi}_{i+1,j}^{n,q}-2\widetilde{\phi}_{i,j}^{n,q}+\widetilde{\phi}_{i-1,j}^{n,q}}{\Delta x^{2}}\right)\cdot(\widetilde{\phi}_{i,j}^{n,q}-\widetilde{\phi}_{i,j}^{n,q-1})\nonumber\\
&\,+\varepsilon^{2}\sum_{j=1}^{N_{y}}\left(\frac{\widetilde{\phi}_{N_{x},j}^{n,q}-\widetilde{\phi}_{N_{x}-1,j}^{n,q}}{\Delta x^{2}}\right)
      \cdot(\widetilde{\phi}_{N_{x},j}^{n,q}-\widetilde{\phi}_{N_{x},j}^{n,q-1})
-\varepsilon^{2}\sum_{j=1}^{N_{y}}\left(\frac{\widetilde{\phi}_{2,j}^{n,q}-\widetilde{\phi}_{1,j}^{n,q}}{\Delta x^{2}}\right)
      \cdot(\widetilde{\phi}_{1,j}^{n,q}-\widetilde{\phi}_{1,j}^{n,q-1})\nonumber\\
&\,-\frac{\varepsilon^{2}}{2}\sum_{i=1}^{N_{x}-1}\sum_{j=1}^{N_{y}}\left[\left(\frac{\widetilde{\phi}_{i+1,j}^{n,q}-\widetilde{\phi}_{i,j}^{n,q}}{\Delta x}\right)
       -\left(\frac{\widetilde{\phi}_{i+1,j}^{n,q-1}-\widetilde{\phi}_{i,j}^{n,q-1}}{\Delta x}\right)\right]^{2},
\label{thf5}
\end{align}
\begin{align}
&\,\frac{\varepsilon^{2}}{2}\sum_{i=1}^{N_{x}}\sum_{j=1}^{N_{y}-1}\left[\left(\frac{\widetilde{\phi}_{i,j+1}^{n,q}-\widetilde{\phi}_{i,j}^{n,q}}{\Delta y}\right)^{2}
       -\left(\frac{\widetilde{\phi}_{i,j+1}^{n,q-1}-\widetilde{\phi}_{i,j}^{n,q-1}}{\Delta y}\right)^{2}\right]\nonumber\\
=&\,-\varepsilon^{2}\sum_{i=1}^{N_{x}}\sum_{j=2}^{N_{y}-1}
   \left(\frac{\widetilde{\phi}_{i,j+1}^{n,q}-2\widetilde{\phi}_{i,j}^{n,q}+\widetilde{\phi}_{i,j-1}^{n,q}}{\Delta y^{2}}\right)\cdot(\widetilde{\phi}_{i,j}^{n,q}-\widetilde{\phi}_{i,j}^{n,q-1})\nonumber\\
&\,+\varepsilon^{2}\sum_{i=1}^{N_{x}}\left(\frac{\widetilde{\phi}_{i,N_{y}}^{n,q}-\widetilde{\phi}_{i,N_{y}-1}^{n,q}}{\Delta y^{2}}\right)
     \cdot(\widetilde{\phi}_{i,N_{y}}^{n,q}-\widetilde{\phi}_{i,N_{y}}^{n,q-1})
-\varepsilon^{2}\sum_{i=1}^{N_{x}}\left(\frac{\widetilde{\phi}_{i,2}^{n,q}-\widetilde{\phi}_{i,1}^{n,q}}{\Delta y^{2}}\right)\cdot(\widetilde{\phi}_{i,1}^{n,q}-\widetilde{\phi}_{i,1}^{n,q-1})
\nonumber\\
&\,-\frac{\varepsilon^{2}}{2}\sum_{i=1}^{N_{x}}\sum_{j=1}^{N_{y}-1}\left[\left(\frac{\widetilde{\phi}_{i,j+1}^{n,q}-\widetilde{\phi}_{i,j}^{n,q}}{\Delta y}\right)
       -\left(\frac{\widetilde{\phi}_{i,j+1}^{n,q-1}-\widetilde{\phi}_{i,j}^{n,q-1}}{\Delta y}\right)\right]^{2}.
\label{thf6}
\end{align}
Therefore,
\begin{align}
\frac{\widetilde{\mathcal{E}}^{n,q}-\widetilde{\mathcal{E}}^{n,q-1}}{\Delta x\Delta y}
=&\,\frac{\varepsilon^{2}}{2}\sum_{i=1}^{N_{x}-1}\sum_{j=1}^{N_{y}}\left[\left(\frac{\widetilde{\phi}_{i+1,j}^{n,q}-\widetilde{\phi}_{i,j}^{n,q}}{\Delta x}\right)^{2}
       -\left(\frac{\widetilde{\phi}_{i+1,j}^{n,q-1}-\widetilde{\phi}_{i,j}^{n,q-1}}{\Delta x}\right)^{2}\right]\nonumber\\
&\,+\frac{\varepsilon^{2}}{2}\sum_{i=1}^{N_{x}}\sum_{j=1}^{N_{y}-1}\left[\left(\frac{\widetilde{\phi}_{i,j+1}^{n,q}-\widetilde{\phi}_{i,j}^{n,q}}{\Delta y}\right)^{2}
       -\left(\frac{\widetilde{\phi}_{i,j+1}^{n,q-1}-\widetilde{\phi}_{i,j}^{n,q-1}}{\Delta y}\right)^{2}\right]\nonumber\\
&\,+\sum_{i=1}^{N_{x}}\sum_{j=1}^{N_{y}}\left[F(\widetilde{\phi}_{i,j}^{n,q})-F(\widetilde{\phi}_{i,j}^{n,q-1})\right]\nonumber\\
\le&\,-\varepsilon^{2}\sum_{i=1}^{N_{x}}\sum_{j=1}^{N_{y}}(\Delta\widetilde{\phi})_{i,j}^{n,q}(\widetilde{\phi}_{i,j}^{n,q}-\widetilde{\phi}_{i,j}^{n,q-1})
    +\widetilde{\xi}^{n,q}\sum_{i=1}^{N_{x}}\sum_{j=1}^{N_{y}}F'(\widetilde{\phi}_{i,j}^{n,q})(\widetilde{\phi}_{i,j}^{n,q}-\widetilde{\phi}_{i,j}^{n,q-1})\nonumber\\
=&\,\sum_{i=1}^{N_{x}}\sum_{j=1}^{N_{y}}(\widetilde{\phi}_{i,j}^{n,q}-\widetilde{\phi}_{i,j}^{n,q-1})\widetilde{\mu}_{i,j}^{n,q}
=-\sum_{i=1}^{N_{x}}\frac{\Delta t}{\Delta x}\left(\widetilde{J}_{i+\frac{1}{2},q}^{n,q}-\widetilde{J}_{i-\frac{1}{2},q}^{n,q}\right)\widetilde{\mu}_{i,q}^{n,q}\nonumber\\
=&\,\sum_{i=1}^{N_{x}-1}\frac{\Delta t}{\Delta x}\widetilde{J}_{i+\frac{1}{2},q}^{n,q}\left(\widetilde{\mu}_{i+1,q}^{n,q}-\widetilde{\mu}_{i,q}^{n,q}\right)
=-\Delta t\sum_{i=1}^{N_{x}-1}\widetilde{J}_{i+\frac{1}{2},q}^{n,q}\widetilde{V}_{i+\frac{1}{2},q}^{n,q}\nonumber\\
=&\,-\Delta t\sum_{i=1}^{N_{x}-1}\left[\left(\widetilde{V}_{i+\frac{1}{2},q}^{n,q}\right)^{+}M(\widetilde{\phi}_{i,q}^{n,q},\widetilde{\phi}_{i+1,q}^{n,q})
                        +\left(\widetilde{V}_{i+\frac{1}{2},q}^{n,q}\right)^{-}M(\widetilde{\phi}_{i+1,q}^{n,q},\widetilde{\phi}_{i,q}^{n,q})\right]\widetilde{V}_{i+\frac{1}{2},q}^{n,q}\nonumber\\
\le&\,-\Delta t\sum_{i=1}^{N_{x}-1}\min\left\{M(\widetilde{\phi}_{i,q}^{n,q},\widetilde{\phi}_{i+1,q}^{n,q}),M(\widetilde{\phi}_{i+1,q}^{n,q},\widetilde{\phi}_{i,q}^{n,q})\right\}
     \left|\widetilde{V}_{i+\frac{1}{2},q}^{n,q}\right|^{2}\le 0.\label{thf7}
\end{align}

Similarly, the iterations in each $y$-direction are also dissipative, i.e.
\begin{equation}
\widehat{\mathcal{E}}^{n,p}-\widehat{\mathcal{E}}^{n,p-1}
\le -\Delta t\Delta x\Delta y\sum_{j=1}^{N_{y}-1}\min\left\{M(\widehat{\phi}_{p,j}^{n,p},\widehat{\phi}_{p,j+1}^{n,p}),M(\widehat{\phi}_{p,j+1}^{n,p},\widehat{\phi}_{p,j}^{n,p})\right\}
     \left|\widehat{V}_{p,j+\frac{1}{2}}^{n,p}\right|^{2}\le 0,\label{thf8}
\end{equation}
where
\begin{align}
\widehat{\mathcal{E}}^{n,p}=&\,
\Delta x\Delta y\sum_{i=1}^{N_{x}-1}\sum_{j=1}^{N_{y}}\frac{\varepsilon^{2}}{2}\left(\frac{\widehat{\phi}_{i+1,j}^{n,p}-\widehat{\phi}_{i,j}^{n,p}}{\Delta x}\right)^{2}
+\Delta x\Delta y\sum_{i=1}^{N_{x}}\sum_{j=1}^{N_{y}-1}\frac{\varepsilon^{2}}{2}\left(\frac{\widehat{\phi}_{i,j+1}^{n,p}-\widehat{\phi}_{i,j}^{n,p}}{\Delta y}\right)^{2}\nonumber\\
&+\Delta x\Delta y\sum_{i=1}^{N_{x}}\sum_{j=1}^{N_{y}}F(\widehat{\phi}_{i,j}^{n,p}).\label{thf9}
\end{align}

In summary, we have
\begin{align}
\mathcal{E}^{n+1}-\mathcal{E}^{n}=&\,\widehat{\mathcal{E}}^{n,p}|_{p=_{N_{x}}}-\widetilde{\mathcal{E}}^{n,q}|_{q=0}
=\sum_{p=1}^{N_{x}}(\widehat{\mathcal{E}}^{n,p}-\widehat{\mathcal{E}}^{n,p-1})
   +\sum_{q=1}^{N_{y}}(\widetilde{\mathcal{E}}^{n,q}-\widetilde{\mathcal{E}}^{n,q-1})\nonumber\\
\le&\, -\Delta t\Delta x\Delta y\sum_{p=1}^{N_{x}}\sum_{j=1}^{N_{y}-1}
\min\left\{M(\widehat{\phi}_{p,j}^{n,p},\widehat{\phi}_{p,j+1}^{n,p}),M(\widehat{\phi}_{p,j+1}^{n,p},\widehat{\phi}_{p,j}^{n,p})\right\}\left|\widehat{V}_{p,j+\frac{1}{2}}^{n,p}\right|^{2}\nonumber\\
&\,-\Delta t\Delta x\Delta y\sum_{q=1}^{N_{y}}\sum_{i=1}^{N_{x}-1}
\min\left\{M(\widetilde{\phi}_{i,q}^{n,q},\widetilde{\phi}_{i+1,q}^{n,q}),M(\widetilde{\phi}_{i+1,q}^{n,q},\widetilde{\phi}_{i,q}^{n,q})\right\}\left|\widetilde{V}_{i+\frac{1}{2},q}^{n,q}\right|^{2}\nonumber\\
\leq &\,0.\label{thf10}
\end{align}
\end{proof}

The above theorems shows that the dimensional-splitting technique
can effectively decouple the multi-dimensional problem into a series of one-dimensional discrete problems
while preserving original structural properties. During each time step,
the most expensive computation is to find the inverse of Jacobian matrix needed by some nonlinear equation solvers (e.g. the Newton--Raphson method). For a $d$-dimensional problem, this technique (which can be readily extended to higher dimensional case) may reduce the computational complexity from $\mathcal{O}(N^{d\beta})$ to $\mathcal{O}(dN^{\beta+d-1})$, in which the cost of finding inverse for a $N\times N$ matrix is $\mathcal{O}(N^{\beta})$ for $2<\beta\le 3$~\cite{Strassen69,Coppersmith90}.

\section{Numerical results}
\label{sec4}

In this section, several numerical examples will be presented to verify the structure-preserving property of our scheme, and to check that the dynamic process of degenerate Cahn--Hilliard equation with Flory--Huggins potential at low temperature is surface diffusion. Unless otherwise specified, the time step $\Delta t$, the spatial step $\Delta x$ ($=\Delta y$) and the interface parameter $\varepsilon$ are set as $10^{-4}$, $0.004$ and $0.02$ respectively.
The initial value of the phase variable is chosen be
\begin{equation}
\phi(\bm{x},t)|_{t=0}=\lambda\tanh\left(\frac{\mathrm{dist}(\bm{x},\Gamma)}{\sqrt{2}\varepsilon}\right),\quad\bm{x}\in\Omega,\label{eqz1}
\end{equation}
where $\Gamma$ represents some curve/surface in the domain $\Omega$, $\mathrm{dist}(\bm{x},\Gamma)$ denotes the signed distance from point $\bm{x}$ to $\Gamma$
and $\lambda$ is set as $1-10^{-4}$ to ensure that the absolute value of initial value of the phase variable does not exceed $1$ (i.e., $\|\phi^{0}\|_{\infty}<1$).
The critical temperature $\theta_{c}$ in the logarithmic  Flory--Huggins potential~\eqref{eq2} is chosen to be $1$ and we select $M(\phi)=1-\phi^{2}$ as the phase-dependent mobility function in numerical simulation.
For the nonlinear equations appeared in our scheme, we shall perform a Newton-type iteration to solve them at each time step.

\subsection{Boundedness and energy stability}

Following the example given in~\cite{Chen19}, the first initial condition we consider is the random form defined in $\Omega=[0,1]^{2}$ as
\begin{equation}
\phi(x,y,t)|_{t=0}=0.2+0.05\times\mathrm{Rand}(x,y) \label{eqz2}
\end{equation}
where $\mathrm{Rand}(x,y)$ represents the uniform random distribution in $[-1,1]$.


\begin{figure}[!htp]
\centering
\includegraphics[width=14cm]{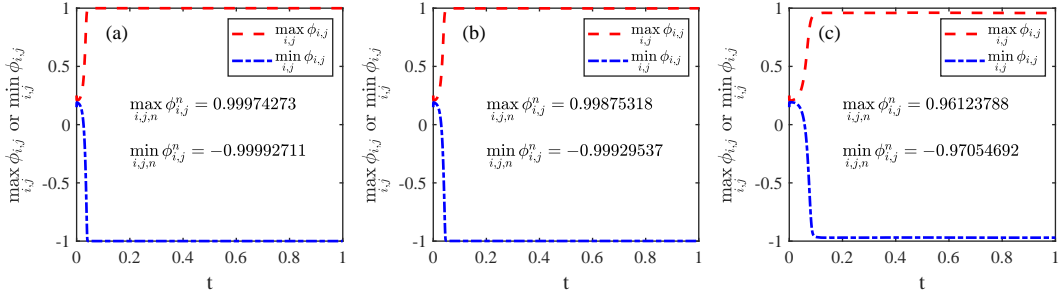}
\caption{The evolution of $\max\limits_{i,j}\phi_{i,j}$ and $\min\limits_{i,j}\phi_{i,j}$ with different $\theta$
         under initial condition~\eqref{eqz2}, where (a) $\theta=0.25$, (b) $\theta=0.30$ and (c) $\theta=0.50$.}\label{frand2}
\end{figure}

\begin{figure}[!htp]
\centering
\includegraphics[width=14cm]{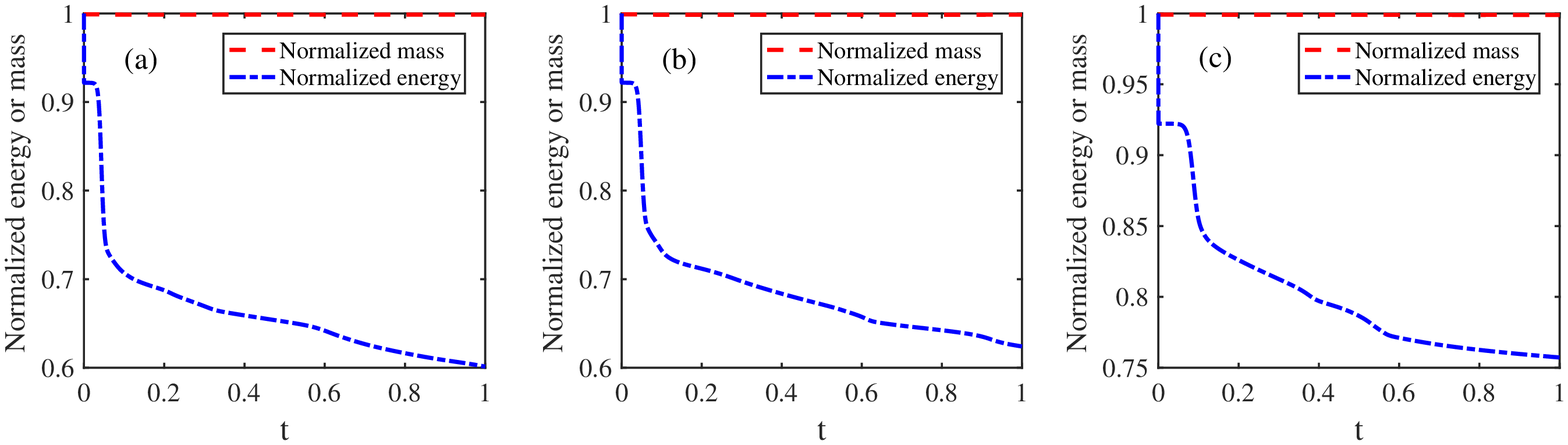}
\caption{The evolution of the normalized energy and mass with different $\theta$ under initial condition~\eqref{eqz2},
         where (a) $\theta=0.25$, (b) $\theta=0.30$ and (c) $\theta=0.50$.}\label{frand3}
\end{figure}

\begin{figure}[!htp]
\centering
\includegraphics[width=14cm]{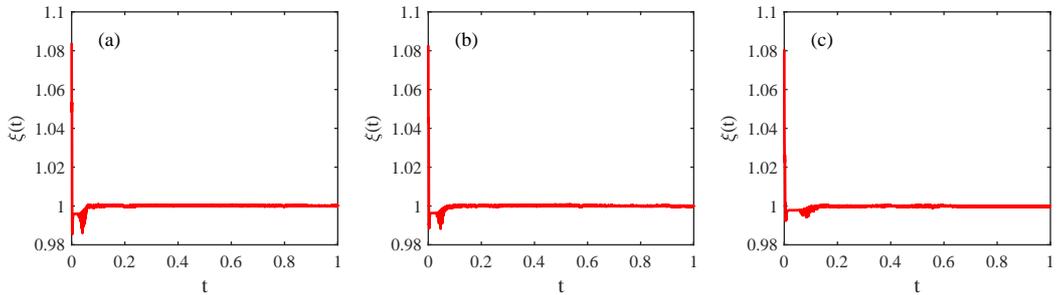}
\caption{The evolution of the Lagrange multiplier $\xi(t)$ with different $\theta$ under initial condition~\eqref{eqz2},
         where (a) $\theta=0.25$, (b) $\theta=0.30$ and (c) $\theta=0.50$.}\label{frand4}
\end{figure}

The boundedness of the proposed scheme is verified in Fig.~\ref{frand2}, in which its maximum (or minimum) approaches $1$ (or $-1$) as $\theta\to 0$. Next, Fig.~\ref{frand3} confirms the principle of mass conservation and energy dissipation.
Finally, we can see from Fig.~\ref{frand4} that the value of Lagrange multiplier $\xi(t)$ is always near $1$ during the evolution,
which is consistent with the theoretical expectation.

To further examine the effect of time step $\Delta t$ on the boundedness and energy stability of numerical solutions,
we select following four-leaved closed curve~\cite{Tang19,Bretin22} at the center of $\Omega=[0,1]^{2}$ as initial condition,
\begin{equation}
\rho=\frac{2+\cos(4\alpha)}{8},\quad\text{where}\quad\rho=\sqrt{x^{2}+y^{2}}\quad\text{and}\quad\alpha=\arctan\frac{y}{x},\label{eqz3}
\end{equation}
and fix $\theta=0.3$ in logarithmic Flory--Huggins potential.


\begin{figure}[htp]
\centering
\includegraphics[width=14cm]{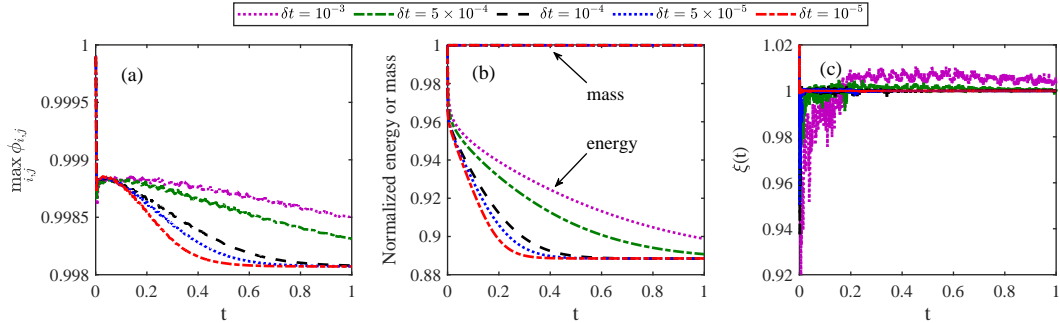}
\caption{The evolution of (a) $\max\limits_{i,j}\phi_{i,j}$, (b) normalized energy and mass and (c) Lagrange multiplier $\xi(t)$
        with different time step size $\Delta t$ under initial condition~\eqref{eqz3}.} \label{fclosed2}
\end{figure}

\begin{table}[!htp]
\renewcommand\arraystretch{1.5}
\centering
\caption{The area change rate $\delta S$ at $t=1$ under initial condition~\eqref{eqz3}.}\label{table1}
\tabcolsep 0.24in
\begin{tabular}[c]
{cccccc} \hline
&$\Delta t=10^{-3}$ &$\Delta t=5\times 10^{-4}$ &$\Delta t=10^{-4}$ &$\Delta t=5\times 10^{-5}$ &$\Delta t=10^{-5}$\\ \cline{2-6}
$\delta S$    &$-0.2996\%$  &$-0.2395\%$  &$-0.1625\%$  &$-0.1388\%$  &$-0.1158\%$\\ \hline
\end{tabular}
\end{table}

The evolution of the maximum/minimum value of phase variable, normalized energy and mass, and Lagrange multiplier $\xi(t)$ are given in Fig.~\ref{fclosed2}, which numerically verifies that the structure-preserving properties of our scheme are unconditionally satisfied.
Finally, the area change rate $\delta S=\frac{S(t)-S(0)}{S(0)}$ at $t=1$ is given in Table~\ref{table1},
it can be seen that a smaller $\Delta t$ leads to a smaller closed area loss,
which indicates that in practice, a small time step may still help us improve the accuracy.

\subsection{Surface diffusion}
This subsection focus on numerically checking that the sharp-interface limit of the Cahn--Hilliard equation with degenerate mobility and Flory--Huggins potential at low temperature  is surface diffusion.

\begin{figure}[htp]
\centering
\includegraphics[width=14cm]{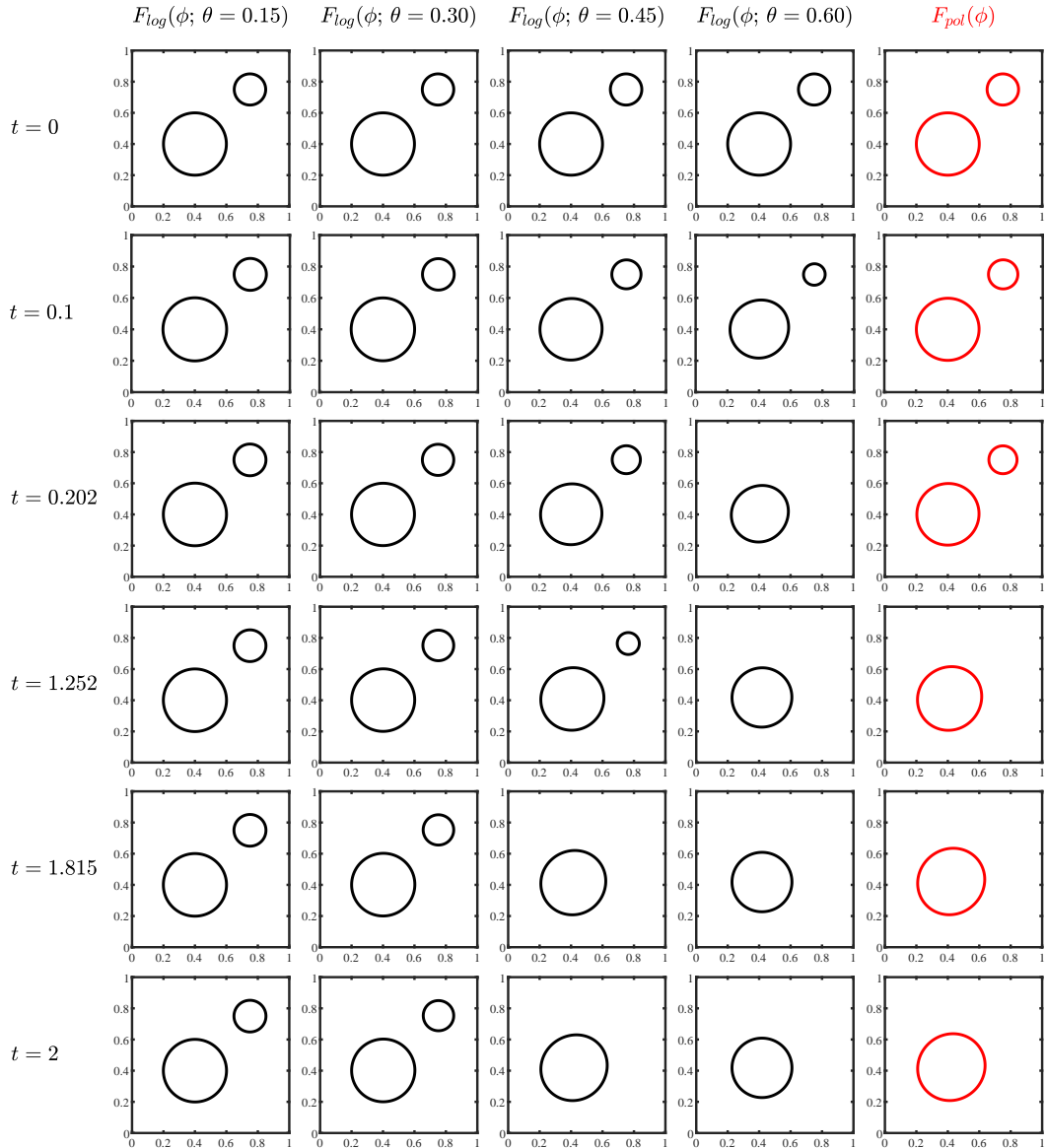}
\caption{Several snapshots of simulating isolated large and small circles with different potential energy densities.}\label{fcircle0}
\end{figure}

\begin{figure}[!htp]
\centering
\includegraphics[width=14cm]{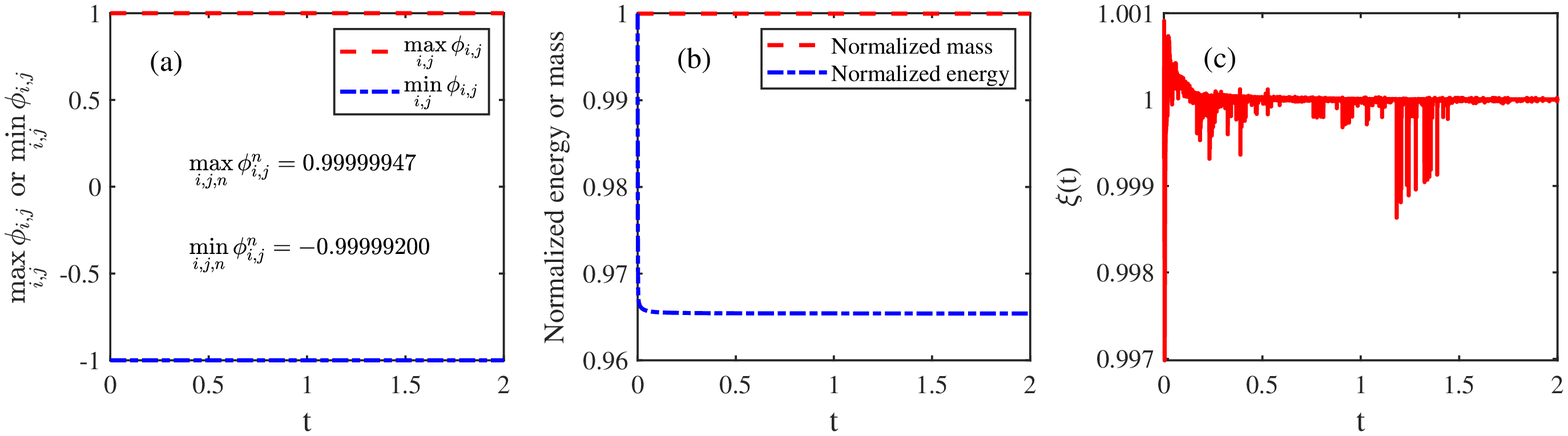}
\caption{The evolution of (a) $\max\limits_{i,j}\phi_{i,j}$ and $\min\limits_{i,j}\phi_{i,j}$, (b) normalized energy and mass and
        (c) Lagrange multiplier $\xi(t)$ under logarithmic potential $F_{log}(\phi)$ at $\theta=0.15$ in Fig.~\ref{fcircle0}.}\label{fcircle1}
\end{figure}

\begin{figure}[!htp]
\centering
\includegraphics[width=14cm]{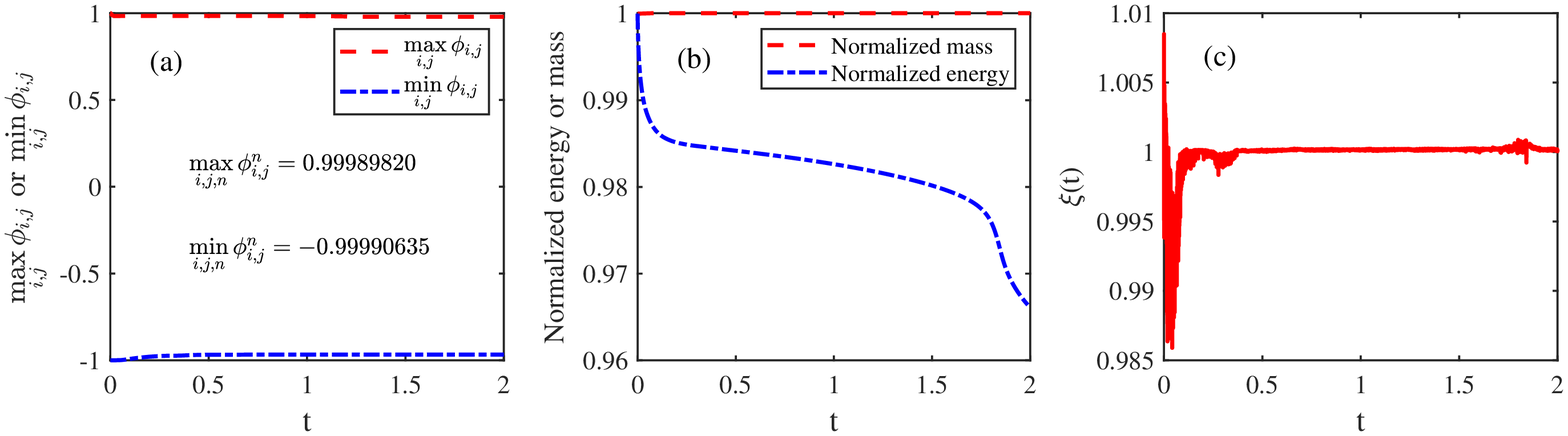}
\caption{The evolution of (a) $\max\limits_{i,j}\phi_{i,j}$ and $\min\limits_{i,j}\phi_{i,j}$, (b) normalized energy and mass and
        (c) Lagrange multiplier $\xi(t)$ under logarithmic potential $F_{log}(\phi)$ at $\theta=0.45$ in Fig.~\ref{fcircle0}.}\label{fcircle3}
\end{figure}

\begin{figure}[!htp]
\centering
\includegraphics[width=14cm]{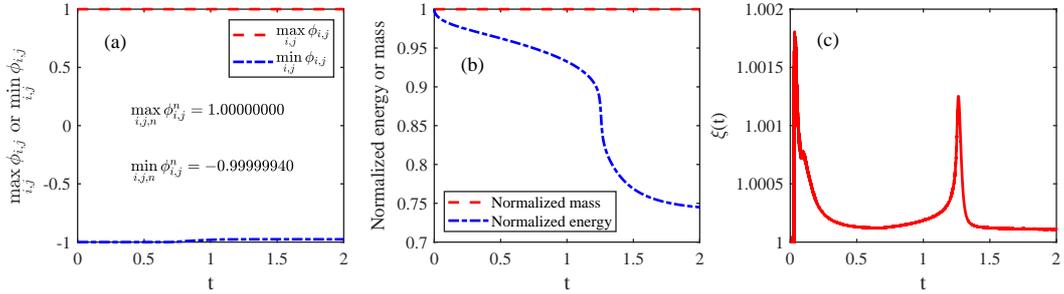}
\caption{The evolution of (a) $\max\limits_{i,j}\phi_{i,j}$ and $\min\limits_{i,j}\phi_{i,j}$, (b) normalized energy and mass and
        (c) Lagrange multiplier $\xi(t)$ under polynomial potential $F_{pol}(\phi)$ in Fig.~\ref{fcircle0}.}\label{fcircle5}
\end{figure}

To begin with, we first consider the evolution of two isolated circles with different sizes under different potential,
one of which has center $(x_{1},y_{1})=(0.4,0.4)$ with radius $r_{1}=0.2$,
and the other one has center $(x_{2},y_{2})=(0.75,0.75)$ with radius $r_{2}=0.1$.
Fig.~\ref{fcircle0} represents several snapshots of solutions under logarithmic potential $F_{log}(\phi)$ with different absolute temperatures $\theta$ and solutions under polynomial potential $F_{pol}(\phi)$.
It can be clearly observed that the process caused by logarithmic potential at large absolute temperatures $\theta$
and polynomial potential violates the main characteristics of surface diffusion,
which is consistent with the theoretical results~\cite{Lee16,Cahn96}.
On the other hand, when the absolute temperature $\theta$ is small enough (e.g., $\theta\le 0.30$ in this simulation),
the kinetic process accords well with the surface diffusion.
Furthermore, the total area change rate $\delta S$ at $t=2$ is also examined. As shown in Table~\ref{tablec}, the area changes very little when $\theta\le 0.30$, which again verifies that a sufficiently low absolute temperature can capture the main feature of surface diffusion.
Finally, the evolution of maximum/minimum value of phase variable, normalized energy and mass, and Lagrange multiplier
under several potential are plotted in Figs.~\ref{fcircle1}-\ref{fcircle5}, respectively,
which are consistent with theoretical expectation.

\begin{table}[!htp]
\renewcommand\arraystretch{1.5}
\centering
\caption{The total area change rate $\delta S$ at $t=2$ under initial two different size circles.}\label{tablec}
\tabcolsep 0.17in
\begin{tabular}[c]
{cccccc} \hline
&$F_{log}(\phi;\,\theta=0.15)$ &$F_{log}(\phi;\,\theta=0.30)$ &$F_{log}(\phi;\,\theta=0.45)$ &$F_{log}(\phi;\,\theta=0.60)$ &$F_{pol}(\phi)$\\ \cline{2-6}
$\delta S$    &$0.1397\%$  &$-0.2829\%$  &$-11.1564\%$  &$-27.2918\%$  &$-8.2058\%$\\ \hline
\end{tabular}
\end{table}

\begin{figure}[htp]
\centering
\includegraphics[width=14cm]{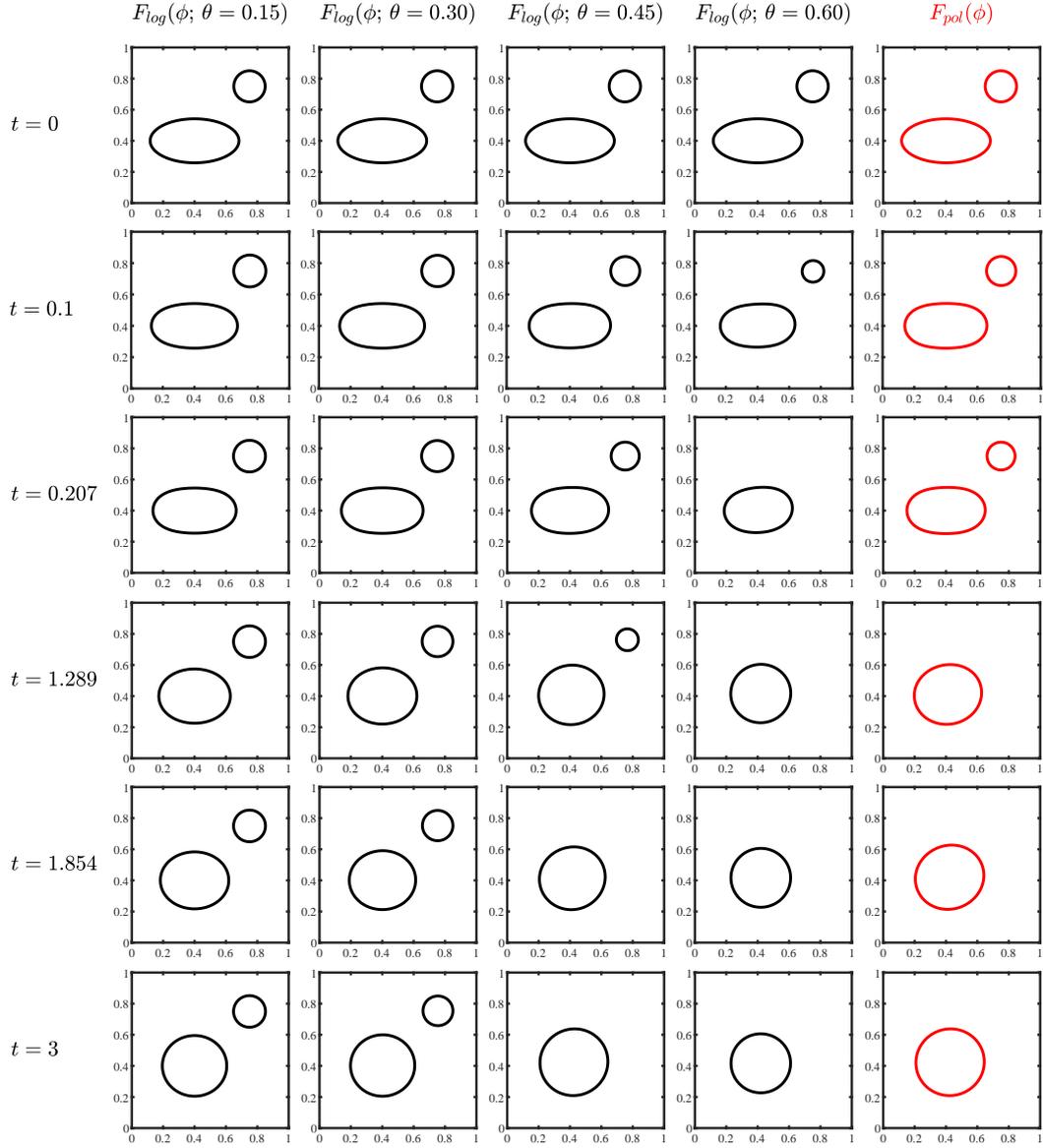}
\caption{Several snapshots of simulating isolated ellipse and small circle with different potential energy densities.}\label{fellipse0}
\end{figure}

\begin{figure}[!htp]
\centering
\includegraphics[width=14cm]{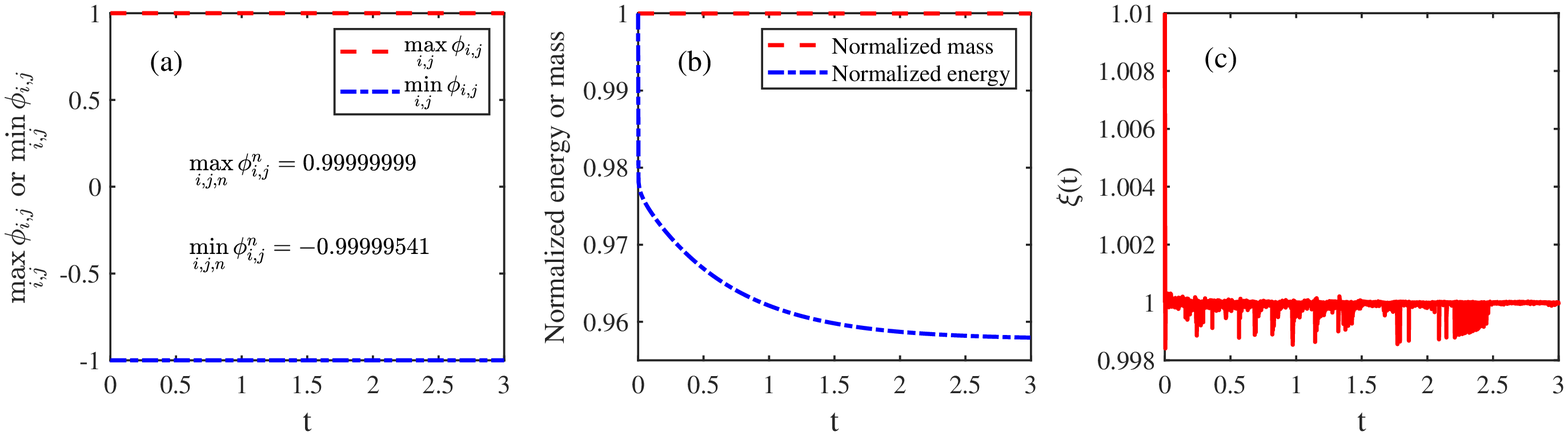}
\caption{The evolution of (a) $\max\limits_{i,j}\phi_{i,j}$ and $\min\limits_{i,j}\phi_{i,j}$, (b) normalized energy and mass and
        (c) Lagrange multiplier $\xi(t)$ under logarithmic potential $F_{log}(\phi)$ at $\theta=0.15$ in Fig.~\ref{fellipse0}.}\label{fellipse1}
\end{figure}

\begin{figure}[!htp]
\centering
\includegraphics[width=14cm]{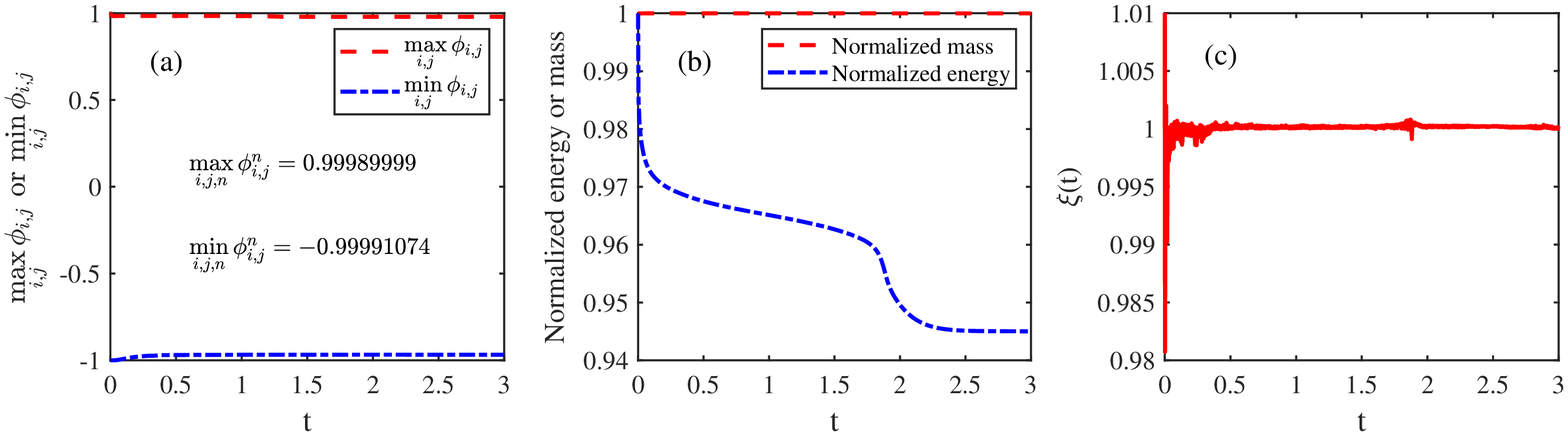}
\caption{The evolution of (a) $\max\limits_{i,j}\phi_{i,j}$ and $\min\limits_{i,j}\phi_{i,j}$, (b) normalized energy and mass and
        (c) Lagrange multiplier $\xi(t)$ under logarithmic potential $F_{log}(\phi)$ at $\theta=0.45$ in Fig.~\ref{fellipse0}.}\label{fellipse3}
\end{figure}

\begin{figure}[!htp]
\centering
\includegraphics[width=14cm]{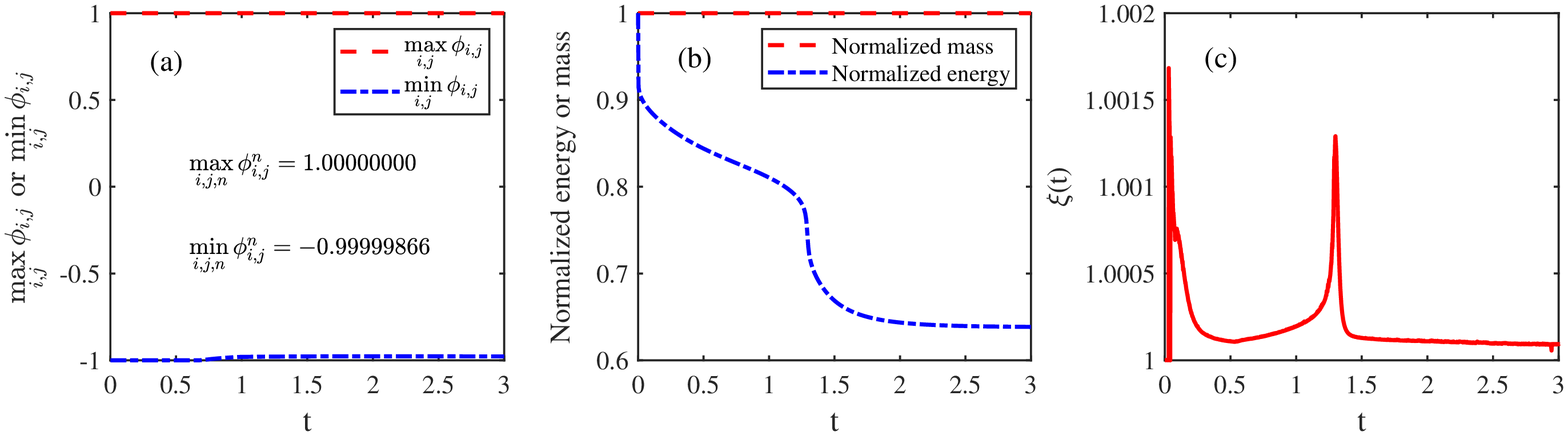}
\caption{The evolution of (a) $\max\limits_{i,j}\phi_{i,j}$ and $\min\limits_{i,j}\phi_{i,j}$, (b) normalized energy and mass and
        (c) Lagrange multiplier $\xi(t)$ under polynomial potential $F_{pol}(\phi)$ in Fig.~\ref{fellipse0}.}\label{fellipse5}
\end{figure}

As the second example, we simply replace the larger circle in the previous example by an ellipse with the major semi axis $r_{a}=\frac{\sqrt{2}}{5}$
and the minor semi axis $r_{b}=\frac{\sqrt{2}}{10}$ (that is, keeping their initial areas equal), leaving the other conditions unchanged,
to test the influence of the ellipse on small circle during the evolution process.
The results of Fig.~\ref{fellipse0} illustrate that when $\theta\le 0.30$, the ellipse gradually evolves into a circle and the small circle remains stable, while in other case the small circle is completely absorbed.
Furthermore, in Table~\ref{tablee}, we quantitatively examine the total area change rate $\delta S$ at $t=3$.
It can be seen that the change of total area is rather small when $\theta\le 0.30$, which again verifies that
the kinetic process of logarithmic potential at low temperature conforms to the geometric characteristics of surface diffusion.
Finally, Figs.~\ref{fellipse1}-\ref{fellipse5} present the evolution of related parameters (e.g., maximum/minimum value of phase variable, normalized energy and mass, Lagrange multiplier) under different potentials, which again numerically illustrate the bound-preserving, mass-conservation and energy-stable properties of the proposed scheme.

\begin{table}[!htp]
\renewcommand\arraystretch{1.5}
\centering
\caption{The total area change rate $\delta S$ at $t=3$ under initial an ellipse and a small circle.}\label{tablee}
\tabcolsep 0.17in
\begin{tabular}[c]
{cccccc} \hline
&$F_{log}(\phi;\,\theta=0.15)$ &$F_{log}(\phi;\,\theta=0.30)$ &$F_{log}(\phi;\,\theta=0.45)$ &$F_{log}(\phi;\,\theta=0.60)$ &$F_{pol}(\phi)$\\ \cline{2-6}
$\delta S$    &$0.7401\%$  &$-0.9826\%$  &$-7.2635\%$  &$-28.0316\%$  &$-6.6129\%$\\ \hline
\end{tabular}
\end{table}

\begin{figure}[htp]
\centering
\includegraphics[width=14cm]{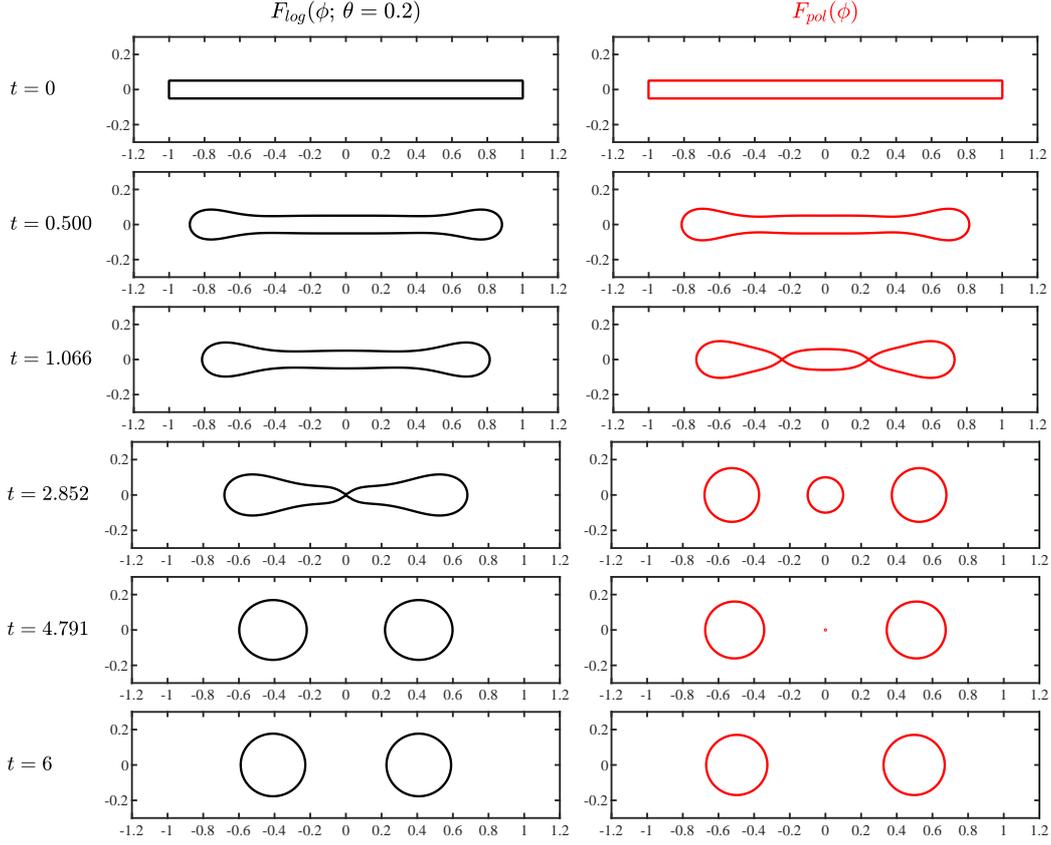}
\caption{Several snapshots of simulating a long rectangle with logarithmic potential $F_{log}(\phi)$ at $\theta=0.2$ (left column)
         and polynomial potential $F_{pol}(\phi)$ (right column).}\label{fpinchoff1}
\end{figure}

\begin{figure}
\centering
\includegraphics[width=14cm]{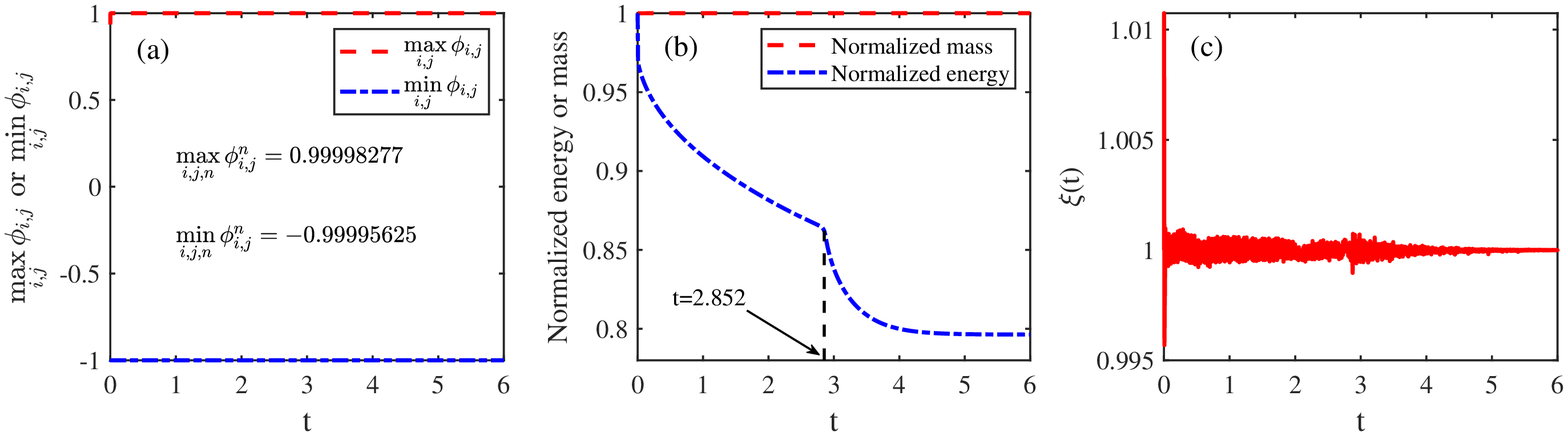}
\caption{The evolution of (a) $\max\limits_{i,j}\phi_{i,j}$ and $\min\limits_{i,j}\phi_{i,j}$, (b) normalized energy and mass and
        (c) Lagrange multiplier $\xi(t)$ under logarithmic potential $F_{log}(\phi)$ at $\theta=0.2$ in Fig.~\ref{fpinchoff1}.}\label{fpinchoff2}
\end{figure}

\begin{figure}
\centering
\includegraphics[width=14cm]{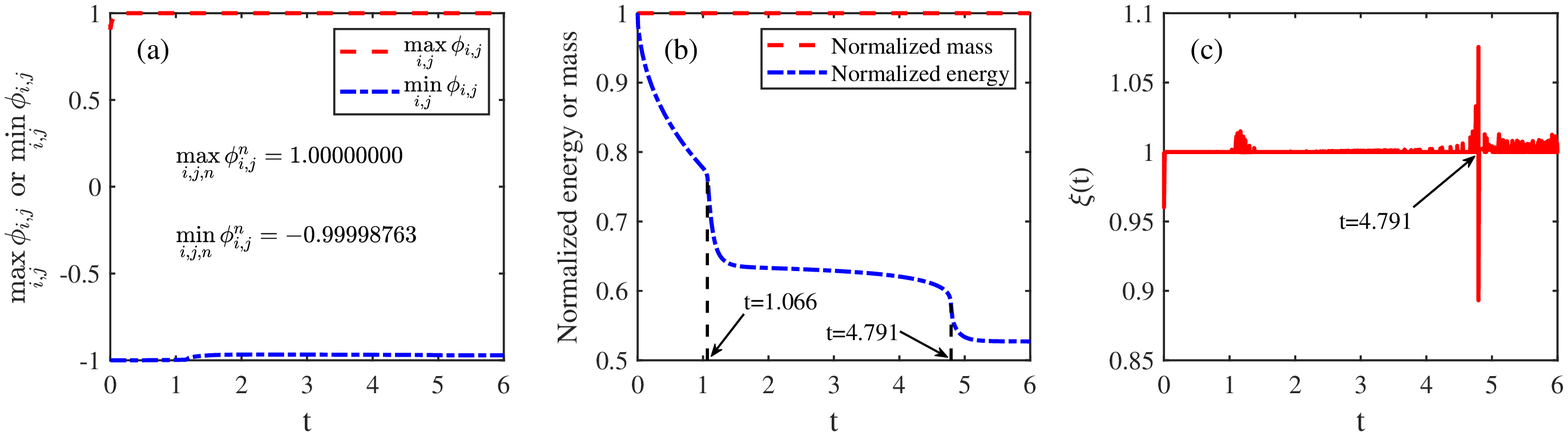}
\caption{The evolution of (a) $\max\limits_{i,j}\phi_{i,j}$ and $\min\limits_{i,j}\phi_{i,j}$, (b) normalized energy and mass and
        (c) Lagrange multiplier $\xi(t)$ under polynomial potential $F_{pol}(\phi)$ in Fig.~\ref{fpinchoff1}.}\label{fpinchoff3}
\end{figure}

At last, since the pinch-off dynamics is also an important phenomenon in the interface evolution problem, we choose to consider a long rectangle with aspect ratio being $20$ as the initial shape. Fig.~\ref{fpinchoff1} depicts the evolution under the logarithmic potential $F_{log}(\phi)$ at $\theta=0.2$ and polynomial potential $F_{pol}(\phi)$.
For polynomial potential $F_{pol}(\phi)$, Fig.~\ref{fpinchoff1} shows that the pinch-off occurs at $t=1.066$
and the long rectangle splits into three independent closed curves.
After that, the middle smaller closed curve is gradually absorbed until it eventually disappears at $t=4.791$.
The total area change rate $\delta S$ equals to $-8.3825\%$ at $t=6$. Obviously, this process is inconsistent with the property of surface diffusion.
As for the case under logarithmic potential $F_{log}(\phi)$ at $\theta=0.2$, the long rectangle splits into two isolated closed curves at $t=2.852$,
and finally evolves into two perfect circles, while the total area change rate $\delta S$ at $t=6$ is only $0.7804\%$.
These results demonstrate that the logarithmic potential at low temperatures with a degenerate mobility can well simulate the pinch-off dynamics of surface diffusion.
In addition, Figs.~\ref{fpinchoff2} and \ref{fpinchoff3} show that the time of pinch-off and disappearance of the middle small closed curve
(only for polynomial potential case) both lead to large changes in the energy curve,
which is consistent with the results of other pinch-off dynamics~\cite{Huang19b,Jiang19}.
Moreover, for polynomial potential, it is interesting that when the middle small closed curve disappears completely,
the Lagrange multiplier will have a relatively bigger deviation from $1$,
which may be caused by the rapid decrease of energy when phase transition occurs,
so the Lagrange multiplier makes a larger correction to the numerical results of fully-implicit scheme
(i.e., the Lagrange multiplier term is removed from the proposed scheme) to ensure the stability of energy.

\section{Conclusion}
\label{seccon}

In this paper, the upwind-scheme and SAV approach are successfully combined to construct
an unconditionally bound-preserving and energy-stable fully-discrete scheme for the Cahn--Hilliard equation with degenerate mobility.
In particular, for a high-dimensional problem, the dimensional-splitting technique is introduced in the SAV framework
for the first time to decouple it into a series of one-dimensional problems while maintaining the original structural properties,
thereby reducing the computational complexity from $\mathcal{O}(N^{d\beta})$ to $\mathcal{O}(dN^{\beta+d-1})$,
which greatly saves computational costs.
Numerical results confirm the boundedness and energy stability of the proposed scheme
for solving degenerate Cahn--Hilliard equations, and successfully capture the main features of surface diffusion numerically
when the absolute temperature in the logarithmic Flory--Huggins potential is sufficiently low.
In future work, we plan to generalize the upwind-SAV approach to other gradient flows with degenerate term,
and the well-posedness and error estimates of this scheme will also be studied.

\begin{acknowledgement}
The numerical calculations in this paper have been done on the supercomputing system in the Supercomputing Center of Wuhan University.
\end{acknowledgement}

%

\section*{Declarations}
{\textbf{Conflict of interest}} ~The authors declare no competing interests.



\end{document}